\documentclass[10pt]{article}
\usepackage{graphicx}
\usepackage{amsmath}
\usepackage[T1]{fontenc}
\usepackage{cmap}
\usepackage{amsfonts}
\usepackage{amssymb}

\newtheorem{theorem}{Theorem}[section]

\newtheorem{corollary}[theorem]{Corollary}

\newtheorem{definition}[theorem]{Definition}

\newtheorem{lemma}[theorem]{Lemma}

\newtheorem{proposition}[theorem]{Proposition}
\newtheorem{remark}[theorem]{Remark}

\newenvironment{proof}[1][Proof]{\textbf{#1.} }{\ \rule{0.5em}{0.5em}}

\begin{document}

\title{Unif{}ication and projectivity in \\Fregean varieties}
\author{Katarzyna S\l omczy\'{n}ska\\ \\Institute of Mathematics, Pedagogical University, \\ul. Podchor\k{a}\.{z}ych 2, 30-084 Krak\'{o}w, Poland \\e-mail: kslomcz@ap.krakow.pl}
\maketitle
\begin{abstract}
In some varieties of algebras one can reduce the question of finding most
general unifiers (mgus) to the problem of the existence of unifiers that
fulfill the additional condition called projectivity. In this paper we study
this problem for Fregean ($1$-regular and orderable) varieties that arise from
the algebraization of fragments of intuitionistic or intermediate logics. We
investigate properties of Fregean varieties, guaranteeing either for a given
unifiable term or for all unifiable terms, that projective unifiers exist. We
indicate the identities which fully characterize congruence permutable Fregean
varieties having projective unifiers. In particular, we show that for such a
variety there exists the largest subvariety that have projective unifiers.\newline

Key words: unification, projectivity, Fregean varieties, intuitionistic logic,
equivalential algebras
\end{abstract}

\section{Introduction}

Equational unification for a variety of algebras is the problem of solving
equations of terms by finding substitutions into their variables, called
unifiers, that make terms equivalent with respect to the variety, see e.g.
\cite{BaaSie94,BaaSny01}. In particular, the problem of the existence of a
most general unifier (an mgu), i.e. such that any other unifier is an instance
of it, is one of the central issues of the unification theory. For certain
varieties the question of searching mgus can be reduced to the problem of
looking for unifiers that obey additional restrictions, namely, the so called
projective unifiers introduced by Ghilardi \cite{Ghi99} (see also
\cite{Wro05}), but in fact known earlier in the unification theory under the
name of reproductive solutions, see e.g. \cite{BaaNip98}. If every unifiable
term has an mgu, we called a variety unitary. In this paper we investigate the
special class of unitary varieties, where every unifiable pair of terms has a
projective unifier and we call them varieties with projective unifiers.

The equational unification has been studied for many varieties related to
logic, see e.g. \cite{Dzi08,Ghi97,Ghi99,Ghi04,Wro95,Wro05,Wro08}. In the
present paper we analyze the problem of unification, and, in particular, the
question of the existence of projective unifiers, for a wide class of such
structures - the congruence permutable (CP) Fregean varieties. The name
Fregean comes from the Frege's idea that sentences with the same logical
values have the same denotation. This idea was formalized by Suszko in
\cite{Suz75} and was an inspiration for Pigozzi \cite{Pig88}, who transferred
the distinction between Fregean and non-Fregean to the field of universal
algebra, see also \cite{CzePig04}. Namely, a variety $\mathcal{V}$ of algebras
with a distinguished constant $1$ is called Fregean if it is $1$-regular, i.e.
congruences of its algebras are uniquely determined by their $1$-cosets, and
congruence orderable, i.e. for every $\mathbf{A}\in\mathcal{V} $ the relation
defined on $A$ by putting $a\leq b$ iff $\Theta_{\mathbf{A}}\left(
1,b\right)  \subset\Theta_{\mathbf{A}}\left(  1,a\right)  $ for $a,b\in A$ is
a partial order.

Many natural examples of Fregean varieties come from the algebraization of
fragments of classical, intuitionistic, or intermediate logics: Boolean
algebras ($\mathbf{CPC}$), Heyting algebras and all its subvarieties
($\mathbf{IPC}$, intermediate logics), Brouwerian semilattices ($\mathbf{IPC}%
$, $\rightarrow$, $\wedge$), equivalential algebras ($\mathbf{IPC}$,
$\leftrightarrow$), Boolean groups ($\mathbf{CPC}$, $\leftrightarrow$) or
Hilbert algebras ($\mathbf{IPC}$, $\rightarrow$). Fregean varieties being
\linebreak $1$-regular are congruence modular, but not necessarily congruence
distributive, which makes studying these classes of algebras more difficult.
On the other hand, the structure of CP Fregean varieties is quite well
understood. In particular, it was proved in \cite{IdzSloWro09} that every CP
Fregean variety consists of algebras that are expansions of equivalential
algebras, i.e. algebras that form an algebraization of the purely
equivalential fragment of the intuitionistic propositional logic.

Subtractive Fregean varieties, that is a Fregean varieties endowed with a
special binary term $\mathsf{s}$ fulfilling $\mathsf{s}(x,x)\approx1$ and
$\mathsf{s}(1,x)\approx x$, form a larger class than CP Fregean varieties,
including e.g. Hilbert algebras. It is not difficult to show that for such
varieties with projective unifiers the unification problems for many equations
$\left\{  s_{i}=t_{i}:i=1,\ldots,k\right\}  $ become equivalent to the
matching problem for single equations of the form $\left\{  p=1\right\}  $
(Proposition \ref{matching}). In other words we are looking for a substitution
that makes $p$ equivalent to $1$ in a given variety. Accordingly, we will
restrict our attention here to such problems and we will say about unifiable
or projective terms if such a substitution exists or is projective,
respectively. We will consider the elementary unification, which means that
the terms in question contain only symbols of the signature of the variety.

Our main results are the following. For a CP Fregean variety we prove a
sufficient condition for a unifiable term to be projective (Theorem
\ref{123}). Using this result we give a twofold characterization of CP Fregean
varieties with projective unifiers (Theorem \ref{1234}). Firstly, we show that
this class can be described by a set of identities. Secondly, we characterize
varieties with projective unifiers by properties of their subdirectly
irreducible algebras. As a consequence, we deduce that for each CP\ Fregean
variety there exists the largest subvariety that has projective unifiers
(e.g., for the variety of Heyting algebras, it is the subvariety of
G\"{o}del-Dummett algebras). Moreover, under the additional assumption of
finite signature, we can characterize all unifiable terms in the variety with
projective unifiers (Proposition \ref{unif}).

The method presented here gives us new proofs of unitarity for several known
cases, as well as it leads to some new results. For instance, it follows from
Theorem \ref{1234} that the variety of equivalential algebras with $0$, being
the algebraic counterpart of the equivalence-negation fragment of the
intuitionistic propositional calculus ($\mathbf{IPC},\leftrightarrow,\lnot$),
is unitary. All these results can be partially generalized to subtractive
Fregean varieties (Proposition \ref{subpro}, Theorem \ref{subthe}) and so
cover such different algebraic structures as Brouwerian semilattices, Hilbert
algebras or equivalential algebras. However, while in the first two cases one
can use a term defining principal filters to get a projective unifier, in the
general case it is impossible: for example, there is no non-trivial subvariety
of equivalential algebras which has equationally definable principle filters
and the only one with definable principle filters is the variety of Boolean
groups \cite{IdzWro97}.

\section{Unification and projectivity}

Let $\mathcal{V}$ be a non-trivial variety of algebras of signature
$\mathcal{F}$. Let us consider the term algebra $T_{\mathcal{F}}(n)$, where
$n\in\mathbb{N}$ is the cardinality of a finite set of variables $\left\{
x_{1},\ldots,x_{n}\right\}  $, and the quotient term algebra $\mathbf{F}%
_{n}:=T_{\mathcal{F}}(n)/\approx_{\mathcal{V}}$, which is a \textsl{free
algebra} in $\mathcal{V}$ with the set of free generators $\left\{
\mathbf{x}_{i}:i=1,\ldots,n\right\}  $. Given $t\in T_{\mathcal{F}}(n)$ and
$\mathbf{A}\in\mathcal{V}$, we adopt the convention of writing $t^{\mathbf{A}%
}:A^{n}\rightarrow A$ for $n$-ary operation on $A$ determined by $t$ and using
the bold symbol $\mathbf{t}$ to denote the corresponding element
$t/\approx_{\mathcal{V}}\ =t^{\mathbf{F}_{n}}\left(  \mathbf{x}_{1}%
,\ldots,\mathbf{x}_{n}\right)  \in\mathbf{F}_{n}$ if no confusion arises. In
particular, we write $\mathbf{c}:=c^{\mathbf{F}_{n}}$ for $c\in T_{\mathcal{F}%
}(0)$.

If $\mathbf{A}\in\mathcal{V}$ and $a,b\in A$, we denote by $\Theta
_{\mathbf{A}}\left(  a,b\right)  $ the congruence generated by $\left(
a,b\right)  $ in $\mathbf{A}$. When there is no ambiguity we drop the
dependence on $\mathbf{A} $ and write just $\Theta\left(  a,b\right)  $.

\begin{definition}
An algebra $\mathbf{A}$ from a variety $\mathcal{V}$ is called
\textsl{projective} in $\mathcal{V}$ if for every algebras $\mathbf{B}%
,\mathbf{C}\in\mathcal{V}$, any epimorphism (onto homomorphism) $\beta
:\mathbf{B}\rightarrow\mathbf{C}$, and any homomorphism $\gamma:\mathbf{A}%
\rightarrow\mathbf{C}$, there exists a homomorphism $\alpha:\mathbf{A}%
\rightarrow\mathbf{B}$ such that $\gamma=\beta\circ\alpha$.
\end{definition}

It is well known that an algebra is projective iff it is a retract of some
free algebra. Moreover, we have the following simple characterization of
projective quotient algebras of a projective algebra.

\begin{proposition}
\label{projquotient}Let $\mathbf{A}$ be a projective algebra in $\mathcal{V}$
and $\varphi\in\operatorname*{Con}\left(  \mathbf{A}\right)  $. The following
conditions are equivalent:

\begin{enumerate}
\item [(1)]$\mathbf{A}/\varphi$ is a projective algebra in $\mathcal{V}$;

\item[(2)] there exists an endomorphism $\tau:A\rightarrow A$ such that
$\varphi=\ker\tau$ and $\tau^{2}=\tau$;

\item[(3)] there exists an endomorphism $\tau:A\rightarrow A$ such that
$\varphi\subset\ker\tau$ and \linebreak $\tau\left(  a\right)  \equiv
_{\varphi}a$ for every $a\in A$.
\end{enumerate}
\end{proposition}

\begin{proof}
$\left(  1\right)  \Rightarrow\left(  2\right)  $. Let $\pi_{\varphi}$ be a
natural epimorphism from $\mathbf{A}$ to $\mathbf{A}/\varphi$. Then there
exists a section $s:\mathbf{A}/\varphi\rightarrow\mathbf{A}$ such that
$\pi_{\varphi}\circ s=\operatorname*{id}_{\mathbf{A}/\varphi}$. Hence
$\tau:=s\circ\pi_{\varphi}$ is a required homomorphism.

$\left(  2\right)  \Rightarrow\left(  3\right)  $. Let $a\in A$. Then
$\tau\left(  \tau\left(  a\right)  \right)  =\tau\left(  a\right)  $, so
$\left(  \tau\left(  a\right)  ,a\right)  \in\ker\tau=\varphi$, as desired.

$\left(  3\right)  \Rightarrow\left(  1\right)  $. For $\mathbf{B}%
,\mathbf{C}\in\mathcal{V}$, an epimorphism $\beta:\mathbf{B}\rightarrow
\mathbf{C}$ and a homomorphism $\gamma:\mathbf{A}/\varphi\rightarrow
\mathbf{C}$, the projectivity of $\mathbf{A}$ implies that there exists a
homomorphism $\alpha:\mathbf{A}\rightarrow\mathbf{B}$ such that $\gamma
\circ\pi_{\varphi}=\beta\circ\alpha$. Put $\overline{\alpha}\left(
a/\varphi\right)  :=\alpha\left(  \tau\left(  a\right)  \right)  $. Then
$\overline{\alpha}:\mathbf{A}/\varphi\rightarrow\mathbf{B}$ is well-defined
and $\beta\left(  \overline{\alpha}\left(  a/\varphi\right)  \right)
=\beta\left(  \alpha\left(  \tau\left(  a\right)  \right)  \right)
=\gamma\left(  \pi_{\varphi}\left(  \tau\left(  a\right)  \right)  \right)
=\gamma\left(  \pi_{\varphi}\left(  a\right)  \right)  =\gamma\left(
a/\varphi\right)  $ for $a\in A$, as required.\smallskip
\end{proof}

There are various ways to introduce unifiers for the problem $\left\{
s=t\right\}  $ in $\mathcal{V}$. We can treat them either as substitutions
from the set of variables to the term algebra, or to define them as follows:

\begin{definition}
Let $s,t\in T_{\mathcal{F}}(n)$ and $\sigma\in\operatorname*{Hom}\left(
\mathbf{F}_{n},\mathbf{F}_{m}\right)  $, $n,m\in\mathbb{N}$. We say that
$\sigma$ is a $\mathcal{V}$\textsl{-unifier\ for }$\left(  s,t\right)  $
\textit{iff }$\sigma\left(  \mathbf{s}\right)  =\sigma\left(  \mathbf{t}%
\right)  $. We say that $\left(  s,t\right)  $ is $\mathcal{V}$%
\textsl{-unifiable} iff there exists a $\mathcal{V}$-unifier for $\left(
s,t\right)  $.
\end{definition}

Studying the unification problem for Heyting algebras, Ghilardi \cite{Ghi99}
considered a special class of unifiers with `good' properties and called them projective.

\begin{definition}
For $s,t\in T_{\mathcal{F}}(n)$, $\tau\in\operatorname*{Hom}\left(
\mathbf{F}_{n},\mathbf{F}_{n}\right)  $ we say that $\tau$ is a $\mathcal{V}%
$\textsl{-projective unifier\ for }$\left(  s,t\right)  $ iff $\tau\left(
\mathbf{s}\right)  =\tau\left(  \mathbf{t}\right)  $ and $\tau\left(
\mathbf{x}_{i}\right)  \equiv_{{\Theta(}\mathbf{s,t}{)}}\mathbf{x}_{i}$ for
every $i=1,\ldots,n$. We say that $\left(  s,t\right)  $ is $\mathcal{V}%
$\textsl{-projective} iff there exists a $\mathcal{V}$-projective unifier for
$\left(  s,t\right)  $.\footnote{Clearly, if $s,t\in T_{\mathcal{F}}(n)$, then
$s,t\in T_{\mathcal{F}}(k)$ for $k\geq n$. However, it is easy to show that
the definitions of $\mathcal{V}$-unifiability and $\mathcal{V}$-projectivity
do not depend on the choice of $k$.}
\end{definition}

Using Proposition \ref{projquotient} we find a simple relation between
projective problems and projective algebras.

\begin{proposition}
\label{projalg}If $s,t\in T_{\mathcal{F}}(n)$, then $\left(  s,t\right)  $ is
$\mathcal{V}$-projective iff $\mathbf{F}_{n}/{\Theta(}\mathbf{s,t}{)}$ is a
projective algebra.
\end{proposition}

For a substitution $\psi:\left\{  x_{1},\ldots,x_{n}\right\}  \rightarrow
T_{\mathcal{F}}(m)$, $m\in\mathbb{N}$, which one can identify with the
homomorphism from $T_{\mathcal{F}}(n)$ to $T_{\mathcal{F}}(m)$, we denote by
$\boldsymbol\psi$ the corresponding homomorphism from $\mathbf{F}_{n}$ to
$\mathbf{F}_{m}$ defined by $\boldsymbol\psi\left(  \mathbf{x}_{\mathbf{i}%
}\right)  :=\mathbf{p}_{\mathbf{i}}$ if $\psi\left(  x_{i}\right)  =p_{i}$ for
$i=1,\ldots,n$. It is not difficult to show that $\boldsymbol\psi$ is a
$\mathcal{V}$-projective unifier if and only if $\psi$ is a
\textsl{reproductive solution }in the sense of Boolean unification, see
\cite{BaaNip98}, or a \textsl{transparent unifier} introduced by Wro\'{n}ski
\cite{Wro05}:

\begin{proposition}
Let $s,t\in T_{\mathcal{F}}(n)$ and $\psi:\left\{  x_{1},\ldots,x_{n}\right\}
\rightarrow T_{\mathcal{F}}(n)$ be a substitution. Then the following three
conditions are equivalent:

\begin{enumerate}
\item [(1)]$\boldsymbol\psi:\mathbf{F}_{n}\rightarrow\mathbf{F}_{n}$ is a
$\mathcal{V}$-projective unifier for $\left(  s,t\right)  $;

\item[(2)] for every $\mathbf{A}\in\mathcal{V}$, $\underline{a}=(a_{1}%
,\ldots,a_{n})\in A^{n}$

\begin{enumerate}
\item $s^{\mathbf{A}}\left(  \psi\left(  \underline{a}\right)  \right)
=t^{\mathbf{A}}\left(  \psi\left(  \underline{a}\right)  \right)  $;

\item $s^{\mathbf{A}}\left(  \underline{a}\right)  =t^{\mathbf{A}}\left(
\underline{a}\right)  $ implies $\psi\left(  \underline{a}\right)
=\underline{a}$,
\end{enumerate}
\end{enumerate}

where $\psi\left(  \underline{a}\right)  \in A^{n}$ is given by $\left(
\psi\left(  \underline{a}\right)  \right)  _{i}:=\left(  \psi\left(
x_{i}\right)  \right)  ^{\mathbf{A}}\left(  \underline{a}\right)  $ for
$i=1,\ldots,n$;

\begin{enumerate}
\item [(3)]$\mathcal{V}\vDash\psi\left(  s\right)  \approx\psi\left(
t\right)  $,

$\mathcal{V}\vDash s\approx t\rightarrow\bigwedge\nolimits_{i=1}^{n}\left(
\psi\left(  x_{i}\right)  \approx x_{i}\right)  $.
\end{enumerate}
\end{proposition}

\begin{proof}
$\left(  1\right)  \Rightarrow\left(  2\right)  $. Consider the homomorphism
$\pi:\mathbf{F}_{n}\rightarrow\mathbf{A}$ given by $\pi\left(  \mathbf{x}%
_{i}\right)  :=a_{i}$. Then $\left(  \psi\left(  \underline{a}\right)
\right)  _{i}=\left(  \psi\left(  x_{i}\right)  \right)  ^{\mathbf{A}}\left(
\pi\left(  \mathbf{x}_{1}\right)  ,\ldots,\pi\left(  \mathbf{x}_{n}\right)
\right)  =\pi\left(  \boldsymbol\psi\left(  \mathbf{x}_{i}\right)  \right)  $
for $i=1,\ldots,n$. Hence we get $s^{\mathbf{A}}\left(  \psi\left(
\underline{a}\right)  \right)  =s^{\mathbf{A}}\left(  \pi\left(
\boldsymbol\psi\left(  \mathbf{x}_{1}\right)  \right)  ,\ldots,\pi\left(
\boldsymbol\psi\left(  \mathbf{x}_{n}\right)  \right)  \right)  =\pi\left(
s^{\mathbf{F}_{n}}\left(  \boldsymbol\psi\left(  \mathbf{x}_{1}\right)
,\ldots,\boldsymbol\psi\left(  \mathbf{x}_{n}\right)  \right)  \right)
=\pi\left(  \boldsymbol\psi\left(  \mathbf{s}\right)  \right)  $, and
analogously $t^{\mathbf{A}}\left(  \psi\left(  \underline{a}\right)  \right)
=\pi\left(  \boldsymbol\psi\left(  \mathbf{t}\right)  \right)  $. Then
$s^{\mathbf{A}}\left(  \psi\left(  \underline{a}\right)  \right)  =\pi\left(
\boldsymbol\psi\left(  \mathbf{s}\right)  \right)  =\pi\left(  \boldsymbol
\psi\left(  \mathbf{t}\right)  \right)  =t^{\mathbf{A}}\left(  \psi\left(
\underline{a}\right)  \right)  $. Assume now that $s^{\mathbf{A}}\left(
\underline{a}\right)  =t^{\mathbf{A}}\left(  \underline{a}\right)  $, which
implies $\pi\left(  \mathbf{s}\right)  =\pi\left(  \mathbf{t}\right)  $. For
$i=1,\ldots,n$ we have $\left(  \boldsymbol\psi\left(  \mathbf{x}_{i}\right)
,\mathbf{x}_{i}\right)  \in{\Theta(}\mathbf{s,t}{)}\subset\ker\pi$, and hence
$a_{i}=\pi\left(  \mathbf{x}_{i}\right)  =\pi\left(  \boldsymbol\psi\left(
\mathbf{x}_{i}\right)  \right)  =\left(  \psi\left(  \underline{a}\right)
\right)  _{i}$.

$\left(  2\right)  \Rightarrow\left(  1\right)  $. Put $\mathbf{A}%
:=\mathbf{F}_{n}$ and $a_{i}:=\mathbf{x}_{i}$ for $i=1,\ldots,n$. Then we get
from (a) $\boldsymbol
\psi\left(  \mathbf{t}\right)  =t^{\mathbf{A}}\left(  \boldsymbol\psi\left(
\mathbf{x}_{1}\right)  ,\ldots,\boldsymbol
\psi\left(  \mathbf{x}_{n}\right)  \right)  =t^{\mathbf{A}}\left(  \psi\left(
\mathbf{x}_{1},\ldots,\mathbf{x}_{n}\right)  \right)  =s^{\mathbf{A}}\left(
\psi\left(  \mathbf{x}_{1},\ldots,\mathbf{x}_{n}\right)  \right)
=s^{\mathbf{A}}\left(  \boldsymbol\psi\left(  \mathbf{x}_{1}\right)
,\ldots,\boldsymbol
\psi\left(  \mathbf{x}_{n}\right)  \right)  =\boldsymbol\psi\left(
\mathbf{s}\right)  $, and so $\boldsymbol\psi$ is a $\mathcal{V}$-unifier for
$\left(  s,t\right)  $. To show that $\boldsymbol\psi$ is projective put, this
time, $\mathbf{A}:=\mathbf{F}_{n}/{\Theta(}\mathbf{s,t}{)}$ and $a_{i}%
:=\mathbf{x}_{i}/{\Theta(}\mathbf{s,t}{)}$ for $i=1,\ldots,n$. Then
$s^{\mathbf{A}}\left(  \underline{a}\right)  =s^{\mathbf{F}_{n}/{\Theta
(}\mathbf{s,t}{)}}\left(  \mathbf{x}_{i}/{\Theta(}\mathbf{s,t}{)}%
,\ldots,\mathbf{x}_{n}/{\Theta(}\mathbf{s,t}{)}\right)  =\mathbf{s}/{\Theta
(}\mathbf{s,t}{)}=\mathbf{t}/{\Theta(}\mathbf{s,t}{)}=t{^{\mathbf{F}%
_{n}/{\Theta(}\mathbf{s,t}{)}}\left(  \mathbf{x}_{i}/{\Theta(}\mathbf{s,t}%
{)},\ldots,\mathbf{x}_{n}/{\Theta(}\mathbf{s,t}{)}\right)  }=t^{\mathbf{A}%
}\left(  \underline{a}\right)  $, and so, using (b), we obtain $\boldsymbol
\psi\left(  \mathbf{x}_{i}\right)  /{\Theta(}\mathbf{s,t}{)}=\left(
\psi\left(  x_{i}\right)  \right)  ^{\mathbf{F}_{n}/{\Theta(}\mathbf{s,t}{)}%
}\left(  \mathbf{x}_{i}/{\Theta(}\mathbf{s,t}{)},\ldots,\mathbf{x}_{n}%
/{\Theta(}\mathbf{s,t}{)}\right)  =$\linebreak $\left(  \psi\left(
\underline{a}\right)  \right)  _{i}=a_{i}=\mathbf{x}_{i}/{\Theta(}%
\mathbf{s,t}{)}$ for $i=1,\ldots,n $, as desired.

$\left(  2\right)  \Leftrightarrow\left(  3\right)  $. Is obvious.\smallskip
\end{proof}

One of the main problems of the unification theory is the existence of most
general unifiers, that is such unifiers that all other unifiers are their
instances. Namely, $\tau\in\operatorname*{Hom}\left(  \mathbf{F}%
_{n},\mathbf{F}_{m}\right)  $, $m\in\mathbb{N}$, is a \textsl{most general
}$\mathcal{V}$-\textsl{unifier} (\textsl{mgu} for short) for $\left(
s,t\right)  $ iff $\tau$ is a unifier for $\left(  s,t\right)  $ and, for
every $\sigma\in\operatorname*{Hom}\left(  \mathbf{F}_{n},\mathbf{F}%
_{k}\right)  $, $k\in\mathbb{N}$, if $\sigma$ is a $\mathcal{V}$-unifier for
$\left(  s,t\right)  $, then there is $\varphi\in\operatorname*{Hom}\left(
\mathbf{F}_{m},\mathbf{F}_{k}\right)  $ such that $\varphi\circ\tau=\sigma$.
The importance of projective unifiers follows also from the fact that they are
special mgus, where $n=m$ and one can put $\varphi=\sigma$, i.e. $\sigma
\circ\tau=\sigma$. (Unifiers fulfilling this condition are called
\textsl{reproductive}, see \cite{BaaNip98}.) Namely, we have the following:

\begin{proposition}
If $s,t\in T_{\mathcal{F}}(n)$ and $\tau\in\operatorname*{Hom}\left(
\mathbf{F}_{n},\mathbf{F}_{n}\right)  $ is a $\mathcal{V}$-projective unifier
for $\left(  s,t\right)  $, then $\tau$ is an mgu for $\left(  s,t\right)  $.
\end{proposition}

\begin{proof}
Let $\sigma\in\operatorname*{Hom}\left(  \mathbf{F}_{n},\mathbf{F}_{k}\right)
$, $k\in\mathbb{N}$ be such that $\sigma\left(  \mathbf{s}\right)
=\sigma\left(  \mathbf{t}\right)  $. It means that $\left(  \mathbf{s}%
,\mathbf{t}\right)  \in\ker\sigma$, and so ${\Theta(}\mathbf{s,t}{)\subset
}\ker\sigma$. Hence $\left(  \tau\left(  \mathbf{x}_{i}\right)  ,\mathbf{x}%
_{i}\right)  \in\ker\sigma$ and, in consequence, $\sigma\left(  \tau\left(
\mathbf{x}_{i}\right)  \right)  =\sigma\left(  \mathbf{x}_{i}\right)  $ for
every $i=1,\ldots,n$, as desired.
\end{proof}

\begin{definition}
If every $\mathcal{V}$-unifiable pair of terms has a $\mathcal{V}$-projective
unifier, we say that $\mathcal{V}$ \textsl{has projective unifiers}. Clearly,
such $\mathcal{V}$ is \textsl{unitary}, i.e. all $\mathcal{V}$-unifiable pairs
of terms have mgus.
\end{definition}

Let us assume that there is a distinguished constant (denoted by $1$) in the
signature of $\mathcal{V}$. Then one can consider matching problems $\left\{
s=1\right\}  $ and introduce accordingly the notion of $\mathcal{V}$-unifiable
and $\mathcal{V}$-projective terms.

\begin{definition}
Let $\mathcal{V}$ be variety with a distinguished constant $1$. We say that a
term $s$ is $\mathcal{V}$\textsl{-unifiable} (resp. $\mathcal{V}%
$\textsl{-projective}) iff $\left(  s,1\right)  $ is $\mathcal{V}$-unifiable
(resp. $\mathcal{V}$-projective).
\end{definition}

We will apply in the sequel the following technical lemma connecting
projectivity of all unifiable terms with projectivity of some quotient algebras.

\begin{lemma}
\label{proj}Let $\mathcal{V}$ be variety with a distinguished constant $1$.
Then the following two conditions are equivalent:

\begin{enumerate}
\item [(1)]every $\mathcal{V}$-unifiable term is $\mathcal{V}$-projective;

\item[(2)] for every $n\in\mathbb{N}$, $\varphi\in\mathsf{Con}\left(
\mathbf{F}_{n}\right)  $ such that:

\begin{enumerate}
\item $\varphi=\bigvee_{j=1}^{k}\Theta\left(  \mathbf{p}_{j},\mathbf{1}%
\right)  $ for some $p_{1},\ldots,p_{k}\in T_{\mathcal{F}}(n)$, and

\item  if $\left(  \mathbf{s},\mathbf{1}\right)  \in\varphi$, then $s$ is
$\mathcal{V}$-unifiable for every $s\in T_{\mathcal{F}}(n)$,
\end{enumerate}

$\mathbf{F}_{n}/\varphi$ is a projective algebra\textit{.}
\end{enumerate}
\end{lemma}

\begin{proof}
[Proof of Lemma]From Proposition \ref{projalg} we get $\left(  2\right)
\Rightarrow\left(  1\right)  $. We shall prove the reverse implication
$\left(  1\right)  \Rightarrow\left(  2\right)  $ by induction on $k$. We may
assume that $k>1$, as the case $k=1$ is obvious. Let $\varphi=\bigvee
_{j=1}^{k}\Theta\left(  \mathbf{p}_{j},\mathbf{1}\right)  $, where
$p_{1},\ldots,p_{k}\in T_{\mathcal{F}}(n)$, fulfills (2b). Put $\psi
:=\bigvee_{j=1}^{k-1}\Theta\left(  \mathbf{p}_{j},\mathbf{1}\right)  $. Then
$\psi$ fulfills (2), and so, by induction, $\mathbf{F}_{n}/\psi$ is
projective. From Proposition \ref{projquotient} we deduce that there exists an
endomorphism $\sigma:\mathbf{F}_{n}\rightarrow\mathbf{F}_{n}$ such that
$\psi\subset\ker\sigma$ and $\sigma\left(  \mathbf{x}_{i}\right)  \equiv
_{\psi}\mathbf{x}_{i}$ for every $i=1,\ldots,n$. Hence we get immediately
$\sigma\left(  \mathbf{p}_{k}\right)  \equiv_{\psi}\mathbf{p}_{k}$, and
consequently $\left(  \sigma\left(  \mathbf{p}_{k}\right)  ,\mathbf{1}\right)
\in\varphi$. From (2b) it follows that $\sigma\left(  \mathbf{p}_{k}\right)  $
is $\mathcal{V}$-unifiable, and so, by (1), there exists $\tau:\mathbf{F}%
_{n}\rightarrow\mathbf{F}_{n}$ such that $\tau\left(  \sigma\left(
\mathbf{p}_{k}\right)  \right)  =\mathbf{1}$ and $\tau\left(  \mathbf{x}%
_{i}\right)  \equiv_{\Theta\left(  \sigma\left(  \mathbf{p}_{k}\right)
,\mathbf{1}\right)  }\mathbf{x}_{i}$ for every $i=1,\ldots,n$. Then
$\varphi=\psi\vee\Theta\left(  \mathbf{p}_{k},\mathbf{1}\right)  \subset
\ker\left(  \tau\circ\sigma\right)  $ and $\tau\left(  \sigma\left(
\mathbf{x}_{i}\right)  \right)  \equiv_{\Theta\left(  \sigma\left(
\mathbf{p}_{k}\right)  ,\mathbf{1}\right)  }\sigma\left(  \mathbf{x}%
_{i}\right)  \equiv_{\psi}\mathbf{x}_{i}$ that implies $\left(  \tau
\circ\sigma\right)  \left(  \mathbf{x}_{i}\right)  \equiv_{\varphi}%
\mathbf{x}_{i}$ for every $i=1,\ldots,n$. Thus, using Proposition
\ref{projquotient} again, we get that $\mathbf{F}_{n}/\varphi$ is projective,
as desired.\smallskip
\end{proof}

\section{Fregean varieties and equivalential algebras}

To make this paper self-contained we provide in this section all the necessary
information about Fregean varieties. For more details, we refer the reader to
\cite{Slo08,IdzSloWro09}.

A variety $\mathcal{V}$ of algebras with a distinguished constant $1$ is
called \textsl{Fregean} if every $\mathbf{A}\in\mathcal{V}$ is:

\begin{itemize}
\item $1$\textsl{-regular}, i.e. $1/\alpha=1/\beta$ implies $\alpha=\beta$ for
all $\alpha,\beta\in\mathsf{Con}(\mathbf{A})$, and

\item \textsl{congruence orderable}, i.e. $\Theta_{\mathbf{A}}\left(
1,a\right)  =\Theta_{\mathbf{A}}\left(  1,b\right)  $ implies $a=b$ for all
$a,b\in A$.
\end{itemize}

Congruence orderability allows us to introduce a natural partial order on the
universe of every $\mathbf{A}\in\mathcal{V}$ by putting
\[
a\leq b\text{ \ \ iff \ \ }{\Theta_{\mathbf{A}}(}1,b{)}\subseteq
{\Theta_{\mathbf{A}}(}1,a{)}\text{ \ \ for \ \ }a,b\in A\text{.}%
\]
Moreover, it follows from $1$-regularity, that the congruences in $\mathbf{A}
$ are uniquely determined by their $1$--cosets traditionally called
\textsl{filters}. The lattice of all filters $\Phi\left(  \mathbf{A}\right)
:=\left\{  1/\alpha:\alpha\in\mathsf{Con}(\mathbf{A})\right\}  $ is isomorphic
to $\mathsf{Con}(\mathbf{A})$. For $\varphi\in\Phi\left(  \mathbf{A}\right)  $
we denote by $\equiv_{\varphi}$ the congruence corresponding to $\varphi$,
and, for $a\in A$, by $\left[  a\right)  :=1/\Theta_{\mathbf{A}}\left(
1,a\right)  =\left\{  b\in A:b\geq a\right\}  $ the filter generated by $a$.

The following characterization of subdirectly irreducible and simple algebras
in Fregean varieties can be easily deduced from \cite[Lemma 2.1]{IdzSloWro09}.
Let $\mathbf{A}$ be an algebra from a Fregean variety $\mathcal{V}$.

\begin{proposition}
\label{FregSubIrr}~

\begin{enumerate}
\item $\mathbf{A}$ is subdirectly irreducible iff there is the largest
non-unit element $\ast$ in $A$. Then, the monolith $\mu\in\Phi\left(
\mathbf{A}\right)  $ has the form $\mu=\left\{  \ast,1\right\}  $ and all
other cosets with respect to $\equiv_{\mu}$ are one element;

\item $\mathbf{A}$ is simple (i.e. $\left|  \mathsf{Con}\left(  \mathbf{A}%
\right)  \right|  =2$) iff $\left|  A\right|  =2$.
\end{enumerate}
\end{proposition}

Let $\mathsf{Fm}\left(  \mathbf{A}\right)  $ denote the set of all completely
meet irreducible filters in $\mathbf{A}$. For each element $\eta$ of
$\mathsf{Fm}\left(  \mathbf{A}\right)  $ there is a unique filter $\eta^{+}%
\in\Phi\left(  \mathbf{A}\right)  $ such that $\alpha\geq\eta^{+}$ whenever
$\alpha>\eta$ for $\alpha\in\Phi\left(  \mathbf{A}\right)  $. The following
result is a simple consequence of Proposition \ref{FregSubIrr} and standard
facts from universal algebra.

\begin{proposition}
\label{monsubirr}~

\begin{enumerate}
\item  If $\eta$ of $\mathsf{Fm}\left(  \mathbf{A}\right)  $, then
$\mathbf{A}/\eta$ is the subdirectly irreducible algebra with the monolith
$\eta^{+}/\eta=\left\{  \ast_{\eta},1/\eta\right\}  $, where $\ast_{\eta
}:=a/\eta$ for any $a\in\eta^{+}\backslash\eta$. Moreover $b/\eta=b/\eta^{+}$
for $b\notin\eta^{+}$;

\item  If $a\in A\backslash\left\{  1\right\}  $, then there exists $\eta$ of
$\mathsf{Fm}\left(  \mathbf{A}\right)  $ such that $a/\eta=\ast_{\eta}$.
\end{enumerate}
\end{proposition}

Many natural examples of Fregean varieties come from the algebraization of
fragments of intuitionistic or intermediate logics. Among them, equivalential
algebras $\mathcal{E}$ play a special role. By an \textsl{equivalential
algebra} we mean a grupoid $\mathbf{A}=\left(  A,\leftrightarrow,1\right)  $
that is the subreduct of a Heyting algebra with the operation $\leftrightarrow
$ given by $x\leftrightarrow y=\left(  x\rightarrow y\right)  \wedge\left(
y\rightarrow x\right)  $. (We adopt further the convention of associating to
the left and ignoring the symbol of equivalence operation.) This notion was
introduced by Kabzi\'{n}ski and Wro\'{n}ski in \cite{KabWro75} as the
algebraic counterpart of the equivalential fragment of intuitionistic
propositional logic. The variety of equivalential algebras is also definable
by the following identities: $xxy=y$; $xyzz=xz(yz)$, $xy(xzz)(xzz)=xy$ and
$xx=1$. Supplementing these axioms by the identity $x=xyy$ (or by the
associativity law) we obtain the smallest non-trivial subvariety of
$\mathcal{E}$, which coincides with the class of Boolean groups and gives the
algebraic semantics for equivalential fragment of classical logic. For a
fuller treatment of equivalential algebras, see \cite{Slo96,Slo08,IdzSloWro09}.

Equivalential algebras form a paradigm for congruence permutable Fregean
varieties, as the following result (see \cite[Theorem 3.8, Corollary
3.9]{IdzSloWro09}) shows:

\begin{theorem}
\label{FregEqv}Let $\mathcal{V}$ be a congruence permutable Fregean variety.
Then there exists a binary term $\leftrightarrow$ such that for every
$\mathbf{A\in}\ \mathcal{V}$:

\begin{enumerate}
\item $\left(  A,\leftrightarrow,1\right)  $ is an equivalential algebra;

\item $\leftrightarrow$ is a principal congruence term of $\mathbf{A}$, i.e.
$a\equiv_{\varphi}b$ iff $ab\in\varphi$ for every $\varphi\in\Phi\left(
\mathbf{A}\right)  $.
\end{enumerate}
\end{theorem}

From Proposition \ref{FregSubIrr}.1 and Theorem \ref{FregEqv} we can easily
deduce the following useful fact:

\begin{proposition}
\label{monolith}Let $\mathbf{A}$ be a subdirectly irreducible algebra from a
congruence permutable Fregean variety $\mathcal{V}$ with the monolith
$\mu=\left\{  \ast,1\right\}  $. Then

\begin{enumerate}
\item $a\ast=a$ for $a\in A\backslash\left\{  \ast,1\right\}  $;

\item $A\backslash\left\{  \ast\right\}  $ is closed under the equivalence operation.
\end{enumerate}
\end{proposition}

\begin{proof}
(1) Let $a\in A\backslash\left\{  \ast,1\right\}  $. Then $\Theta_{\mathbf{A}%
}\left(  1,\ast\right)  \subset\Theta_{\mathbf{A}}\left(  1,a\right)  $, and
so $\Theta_{\mathbf{A}}\left(  \ast,a\right)  \subset\Theta_{\mathbf{A}%
}\left(  1,a\right)  $. Moreover, $\Theta_{\mathbf{A}}\left(  1,\ast\right)
\subset\Theta_{\mathbf{A}}\left(  \ast,a\right)  $, and in consequence
$\Theta_{\mathbf{A}}\left(  1,a\right)  \subset\Theta_{\mathbf{A}}\left(
\ast,a\right)  $. Hence $\Theta_{\mathbf{A}}\left(  1,a\ast\right)
=\Theta_{\mathbf{A}}\left(  \ast,a\right)  =\Theta_{\mathbf{A}}\left(
1,a\right)  $, as desired.

(2) Assume for contradiction that $\ast=ab$ for some $a,b\in A\backslash
\left\{  \ast,1\right\}  $. Then from (1) $1=\ast aa=abaa=ab$, a contradiction.\smallskip
\end{proof}

Applying Theorem \ref{FregEqv} we can define in every algebra $\mathbf{A}$
from a congruence permutable Fregean variety $\mathcal{V}$ a family of purely
equivalential unary polynomials given by the formula:
\[
\chi_{a}(x):=xaa\text{ for }x,a\in A\text{,}%
\]
with the following properties for $a,b\in A$:

\begin{itemize}
\item $\chi_{a}$ is an idempotent equivalential endomorphism,

\item $\chi_{a}\circ\chi_{b}=\chi_{b}\circ\chi_{a}=\chi_{a}\circ\chi_{ab}$.
\end{itemize}

Note that the polynomial $\chi_{a}(x)$ can be expressed in the language of
Heyting algebra as $\left(  x\rightarrow a\right)  \rightarrow x$. For some
congruence permutable Fregean varieties (e.g. equivalential algebras,
Brouwerian semilattices) these polynomials are endomorphisms, whereas for
others (e.g. Heyting algebras) they are not. These distinction is crucial for
our purposes.\smallskip

A variety $\mathcal{V}$ with the signature $\mathcal{F}$ and a distinguished
constant $1$ is called \textsl{subtractive} iff there is a
\textsl{subtractive} term $\mathsf{s}\in T_{\mathcal{F}}(2)$, i.e. a binary
term such that the identities $\mathsf{s}(x,x)\approx1$ and $\mathsf{s}%
(1,x)\approx x$ are fulfilled in $\mathcal{V}$. This notion was introduced in
\cite{GumUrs84}, see also \cite{Urs94}. Clearly, the class of all subtractive
Fregean varieties is larger than the class of congruence permutable Fregean
varieties, where the equivalence $\leftrightarrow$ serves as the subtractive
term. To be more precise, the existence of a subtractive term in $\mathcal{V}$
is equivalent to a special case of congruences permutability, namely, if
$\alpha$ and $\beta$ are congruences of $\mathbf{A}\in\mathcal{V}$, then
$\left(  1,c\right)  \in\alpha\circ\beta$ iff $\left(  1,c\right)  \in
\alpha\vee\beta$ for $c\in A$. This class was studied by Agliano \cite{Agl01},
who proved that a subtractive Fregean variety has equationally definable
principal congruences iff it is term equivalent to a variety of Hilbert
algebras with additional compatible operations. In the sequel, we will need
the following simple property of subtractive Fregean varieties:

\begin{proposition}
\label{subcon}Let $\mathbf{A}$ be an algebra from a subtractive Fregean
variety $\mathcal{V}$ and let $a,b\in A$. Then $\Theta\left(  a,b\right)
=\Theta\left(  \mathsf{s}(a,b),1\right)  \vee\Theta\left(  \mathsf{s}%
(b,a),1\right)  $, where $\mathsf{s}$ is a subtractive term in $\mathcal{V}$.
\end{proposition}

\begin{proof}
Put $\varphi:=\Theta\left(  \mathsf{s}(a,b),1\right)  \vee\Theta\left(
\mathsf{s}(b,a),1\right)  $. Clearly, $\mathsf{s}(a,b),\mathsf{s}%
(b,a)\equiv_{\Theta\left(  a,b\right)  }\mathsf{s}(b,b)=1$, and, in
consequence, $\varphi\subset\Theta\left(  a,b\right)  $. Moreover,
$a=\mathsf{s}(1,a)=\mathsf{s}(\mathsf{s}(\mathsf{s}(b,a),\mathsf{s}%
(b,a)),a)\equiv_{\varphi}\mathsf{s}(\mathsf{s}(1,\mathsf{s}(b,a)),a)\equiv
_{\Theta\left(  b,1\right)  }\mathsf{s}(\mathsf{s}(1,\mathsf{s}(1,a)),a)=1$.
Hence $\Theta\left(  a,1\right)  \subset\Theta\left(  b,1\right)  \vee\varphi
$, and analogously $\Theta\left(  b,1\right)  \subset\Theta\left(  a,1\right)
\vee\varphi$. Thus $\Theta\left(  a,1\right)  \vee\varphi=\Theta\left(
b,1\right)  \vee\varphi$. It means that $\Theta_{\mathbf{A/}\varphi}\left(
a/\varphi,1/\varphi\right)  =\Theta_{\mathbf{A/}\varphi}\left(  b/\varphi
,1/\varphi\right)  $. Since $\mathbf{A/}\varphi$ is congruence orderable, we
get $a/\varphi=b/\varphi$, and so $\Theta\left(  a,b\right)  \subset\varphi$,
as desired.\smallskip
\end{proof}

In the next sections we will study the problem of the existence of projective
unifiers for congruence permutable and subtractive Fregean varieties. It
follows from Theorem \ref{FregEqv} that for congruence permutable Fregean
variety the unification problem in $\mathcal{V}$ for an arbitrary equation of
the form $\left\{  s=t\right\}  $, where $s,t\in T_{\mathcal{F}}(n)$ for some
$n\in\mathbb{N}$ can be easily reduced to the matching problem for the
equation $\left\{  st=1\right\}  $. It is not true for a subtractive Fregean
variety, however the problem of having projective unifiers for single
equations, or even for finite sets of equations, for such a variety can also
be reduced to the matching problems, as the next proposition shows.

\begin{definition}
We say that the set of equations $\left\{  s_{i}=t_{i}:i=1,\ldots,k\right\}
$, where $s_{i},t_{i}\in T_{\mathcal{F}}(n)$ for some $n\in\mathbb{N}$ is:

$\mathcal{V}$\textsl{-unifiable} iff there exists $\sigma\in
\operatorname*{Hom}\left(  \mathbf{F}_{n},\mathbf{F}_{m}\right)  $,
$m\in\mathbb{N}$ such that\textit{\ }$\sigma\left(  \mathbf{s}_{i}\right)
=\sigma\left(  \mathbf{t}_{i}\right)  $ for $i=1,\ldots,k$;

$\mathcal{V}$\textsl{-projective} iff there exists $\tau\in\operatorname*{Hom}%
\left(  \mathbf{F}_{n},\mathbf{F}_{n}\right)  $ such that\textit{\ }%
$\tau\left(  \mathbf{s}_{i}\right)  =\tau\left(  \mathbf{t}_{i}\right)  $ and
$\tau\left(  \mathbf{x}_{i}\right)  \equiv_{\bigvee_{j=1}^{k}{\Theta
(}\mathbf{s}_{i}\mathbf{,t}_{i}{)}}\mathbf{x}_{i}$ for $i=1,\ldots,k$.
\end{definition}

\begin{proposition}
\label{matching}Let $\mathcal{V}$ be a subtractive Fregean variety. Then the
following conditions are equivalent:

\begin{enumerate}
\item [(1)]every $\mathcal{V}$-unifiable term is $\mathcal{V}$-projective;

\item[(2)] $\mathcal{V}$ has projective unifiers;

\item[(3)] every finite and $\mathcal{V}$-unifiable set of equations is
$\mathcal{V}$-projective.
\end{enumerate}
\end{proposition}

\begin{proof}
[Proof of Theorem]It is obvious that $\left(  3\right)  \Rightarrow\left(
2\right)  \Rightarrow\left(  1\right)  $. To show $\left(  1\right)
\Rightarrow\left(  3\right)  $ take $s_{i},t_{i}\in T_{\mathcal{F}}(n)$ for
some $n\in\mathbb{N}$ and $i=1,\ldots,k$ such that $\left\{  s_{i}%
=t_{i}:i=1,\ldots,k\right\}  $ is $\mathcal{V}$-unifiable. Put $\varphi
:=\bigvee_{j=1}^{k}\Theta\left(  \mathbf{s}_{j},\mathbf{t}_{j}\right)  $. It
follows from Proposition \ref{subcon} that $\varphi=\bigvee_{j=1}^{k}%
\Theta\left(  \mathsf{s}(\mathbf{s}_{j},\mathbf{t}_{j}),\mathbf{1}\right)
\vee\bigvee_{j=1}^{k}\Theta\left(  \mathsf{s}(\mathbf{t}_{j},\mathbf{s}%
_{j}),\mathbf{1}\right)  $. Applying Lemma \ref{proj} we deduce that the
finitely presented algebra $\mathbf{F}_{n}/\varphi$ is projective. Now, it
follows from Proposition \ref{projquotient} that $\left\{  s_{i}%
=t_{i}:i=1,\ldots,k\right\}  $ is $\mathcal{V}$-projective, which completes
the proof.\smallskip
\end{proof}

The conclusion of the above Proposition is rather obvious for any arithmetical
(i.e. congruence permutable and congruence distributive) Fregean variety,
because such a variety must be term equivalent to a variety of Brouwerian
semilattices with some compatible operations, see \cite{Pig88} and
\cite[Corollary 4.1]{IdzSloWro09}. It is rather surprising that the result can
be extended to non-arithmetical Fregean varieties, such as, e.g. equivalential algebras.

\section{Congruence permutable Fregean varieties}

We start from the result that is in fact a simple generalization of the
well-known theorem of Diego, who proved that finitely generated Hilbert
algebras are finite \cite{Die65}. This result has been later extended to other
varieties that come from the algebraization of some fragments of
intuitionistic logic, like Brouwerian semilattices or equivalential algebras.

\begin{theorem}
\label{DieGen}Let $\mathcal{V}$ be a Fregean variety with the finite signature
$\mathcal{F}$. If for every $\mathbf{A}$ subdirectly irreducible in
$\mathcal{V}$ and such that $\left|  A\right|  >2$ the set $A\backslash
\left\{  \ast\right\}  $ is a subuniverse of $\mathbf{A}$, then $\mathcal{V}$
is locally finite.
\end{theorem}

\begin{proof}
We will prove by induction on $n$ that $\left|  F_{n}\right|  <\infty$. From
the finiteness of the signature $\mathcal{F}$ it is enough to show that for
every $n\in\mathbb{N}$ there exists a common finite upper bound for the
cardinality of subdirectly irreducible homomorphic images of $\mathbf{F}_{n}$.
Then there exists only a finite number of non-isomorphic subdirectly
irreducible homomorphic images of $\mathbf{F}_{n}$, and so $\mathbf{F}_{n}$ is
finite as a finitely generated member of a variety generated by a finite
number of finite algebras.

Let $n=0$. Observe that if $\mathbf{F}_{0}$ is non-trivial, then every algebra
$\mathbf{A}$ that is a subdirectly irreducible homomorphic image of
$\mathbf{F}_{0}$ has two elements. Otherwise, $A\backslash\left\{
\ast\right\}  $ is a subuniverse of $\mathbf{A}$, which is impossible, since
$\mathbf{F}_{0}$ has no proper subalgebras. Assume now that $\left|
F_{n}\right|  <\infty$. For $\mu\in\mathsf{Fm}\left(  \mathbf{F}_{n+1}\right)
$, either $\left|  \mathbf{F}_{n+1}/\mu\right|  =2$ or $\left|  \mathbf{F}%
_{n+1}/\mu\right|  >2$. In the latter case $\left(  \mathbf{F}_{n+1}%
/\mu\right)  \backslash\left\{  \ast\right\}  $ is a subuniverse of
$\mathbf{F}_{n+1}/\mu$, and consequently, there is $i=1,\ldots,n+1$ such that
$\mathbf{x}_{i}/\mu=\ast$. Hence and from Proposition \ref{monsubirr} we
deduce that $\mathbf{F}_{n+1}/\mu^{+}\simeq\left(  \mathbf{F}_{n+1}%
/\mu\right)  /\left(  \mu^{+}/\mu\right)  $ has at most $n$ generators and
$\left|  \mathbf{F}_{n+1}/\mu\right|  =1+\left|  \mathbf{F}_{n+1}/\mu
^{+}\right|  $. Thus $\left|  \mathbf{F}_{n+1}/\mu\right|  \leq1+\left|
\mathbf{F}_{n}\right|  $, and now the assertion follows from the induction hypothesis.\smallskip
\end{proof}

The next theorem gives two equivalent sufficient conditions for a unifiable
term $t$ in a congruence permutable Fregean variety $\mathcal{V}$ to be
projective. The first has the form of identity in $\mathcal{V}$. This identity
expresses the fact that for every algebra $\mathbf{A}$ in $\mathcal{V}$ and
$a\in A$, the unary polynomial $\chi_{a}:A\rightarrow A$, defined by $\chi
_{a}(x)=axx$, for $x\in A$, preserves the operation $t^{\mathbf{A}}$. The
second says that for every non-simple subdirectly irreducible algebra
$\mathbf{A}$ in $\mathcal{V}$ the universe of $\mathbf{A}$ with the largest
non-unit element removed is closed under the operation $t^{\mathbf{A}}$.

\begin{theorem}
\label{123}Let $\mathcal{V}$ be a congruence permutable Fregean variety and
$t\in T_{\mathcal{F}}(n)$. Then the following conditions are equivalent:

\begin{enumerate}
\item [(1)]$\mathcal{V}\models t\left(  x_{1}yy,\ldots,x_{n}yy\right)  \approx
t\left(  x_{1},\ldots,x_{n}\right)  yy$ ;

\item[(2)] for every subdirectly irreducible $\mathbf{A}$ in $\mathcal{V}$
such that $\left|  A\right|  >2$, and for all $a_{1},\ldots,a_{n}\in A$%
\[
\left(  t^{\mathbf{A}}\left(  a_{1},\ldots,a_{n}\right)  =\ast\right)
\Rightarrow\left(  \exists i=1,\ldots,n:a_{i}=\ast\right)  \text{ .}%
\]
Moreover, if $t$ fulfills (1) or (2), then

\item[(3)] if $t$ is $\mathcal{V}$-unifiable, then $t$ is $\mathcal{V}$-projective.
\end{enumerate}
\end{theorem}

Before we begin the proof of the theorem we need the following lemma:

\begin{lemma}
Let $\mathbf{A}$ be subdirectly irreducible in $\mathcal{V}$ with $\left|
A\right|  >2$. Then for all $a_{1},\ldots,a_{n}\in A\backslash\left\{
\ast\right\}  $ there exists $p\in A\backslash\left\{  \ast,1\right\}  $ such
that $a_{i}pp=a_{i}$ for every $i=1,\ldots,n$.
\end{lemma}

\begin{proof}
[Proof of Lemma]The proof is by induction on $n$. For $n=1$ we put $p=a_{1}$
if $a_{1}\neq1$ and $p=c$, where $c\in A\backslash\left\{  \ast,1\right\}  $
if $a_{1}=1$. Assume that $n>1$. Take $q\in A\backslash\left\{  \ast
,1\right\}  $ such that $a_{i}qq=a_{i}$ for every $i=1,\ldots,n-1$. It is
enough to consider the case $a_{n}qq\neq a_{n}$. Put $p=a_{n}qqa_{n}\neq1$.
From Proposition \ref{monolith}.2 we deduce that $p\neq\ast$. For
$i=1,\ldots,n-1$ we have $a_{i}pp=a_{i}qq\left(  a_{n}qqa_{n}\right)  \left(
a_{n}qqa_{n}\right)  =a_{i}\left(  a_{n}qqa_{n}\right)  \left(  a_{n}%
qqa_{n}\right)  qq=\left(  a_{i}qq\right)  \left(  \left(  a_{n}%
qqa_{n}\right)  qq\right)  \left(  \left(  a_{n}qqa_{n}\right)  qq\right)
=a_{i}qq=a_{i}$. Moreover, $a_{n}pp=a_{n}\left(  a_{n}qqa_{n}\right)  \left(
a_{n}qqa_{n}\right)  =\left(  a_{n}qqa_{n}a_{n}\right)  \left(  a_{n}%
qqa_{n}\right)  =\left(  a_{n}qq\right)  \left(  a_{n}qqa_{n}\right)
=a_{n}\left(  a_{n}qq\right)  \left(  a_{n}qq\right)  =a_{n}$.\medskip
\end{proof}

\begin{proof}
[Proof of Theorem]$\left(  1\right)  \Rightarrow\left(  2\right)  $. On the
contrary, assume that there exists $\mathbf{A}$ subdirectly irreducible in
$\mathcal{V}$, $\left|  A\right|  >2$ and $a_{1},\ldots,a_{n}\in
A\backslash\left\{  \ast\right\}  $ such that $t^{\mathbf{A}}\left(
a_{1},\ldots,a_{n}\right)  =\ast$. From Lemma there exists a $p\in
A\backslash\left\{  \ast,1\right\}  $ such that $a_{i}pp=a_{i}$ for every
$i=1,\ldots,n$. Hence $\ast=t^{\mathbf{A}}\left(  a_{1},\ldots,a_{n}\right)
=t^{\mathbf{A}}\left(  a_{1}pp,\ldots,a_{n}pp\right)  =t^{\mathbf{A}}\left(
a_{1},\ldots,a_{n}\right)  pp$. On the other hand, from Proposition
\ref{monolith}.1 we get \linebreak $t^{\mathbf{A}}\left(  a_{1},\ldots
,a_{n}\right)  pp=\ast pp=1$, a contradiction.

$\left(  2\right)  \Rightarrow\left(  1\right)  $. Assume to the contrary that
there exists $\mathbf{C}\in\mathcal{V}$ such that the identity in $(1)$ is not
fulfilled. Then we find an algebra $\mathbf{A}$ in $\mathcal{V}$ being the
subdirectly irreducible homomorphic image of $\mathbf{C}$\ and $a_{1}%
,\ldots,a_{n},b\in A$ such that $\left(  t^{\mathbf{A}}\left(  a_{1}%
bb,\ldots,a_{n}bb\right)  ,t^{\mathbf{A}}\left(  a_{1},\ldots,a_{n}\right)
bb\right)  $ generates the monolith of $\mathbf{A}$. It follows from
Proposition \ref{FregSubIrr}.1 that $\left\{  t^{\mathbf{A}}\left(
a_{1}bb,\ldots,a_{n}bb\right)  ,t^{\mathbf{A}}\left(  a_{1},\ldots
,a_{n}\right)  bb\right\}  =\left\{  \ast,1\right\}  $. Hence and from
Proposition \ref{monolith}.1 we deduce that $b\notin\left\{  \ast,1\right\}
$, and so $\left|  A\right|  >2$. Moreover, applying again Proposition
\ref{monolith} and the fact that $b\notin\left\{  \ast,1\right\}  $ we get
$ubb\neq\ast$ for every $u\in A$. In consequence, $t^{\mathbf{A}}\left(
a_{1},\ldots,a_{n}\right)  bb=1$, $t^{\mathbf{A}}\left(  a_{1}bb,\ldots
,a_{n}bb\right)  =\ast$ and $a_{i}bb\neq\ast$ for $i=1,\ldots,n$, which
contradicts (2).

$\left(  1\right)  \Rightarrow\left(  3\right)  $. Since $t$ is unifiable, we
can find an endomorphism $f:\mathbf{F}_{n}\rightarrow\mathbf{F}_{n}$ such that
$f\left(  \mathbf{t}\right)  =\mathbf{1}$. We can assume that $n>0$, since
otherwise $\mathbf{t}=1$ and we are done.

We divide $\mathsf{Fm}\left(  \mathbf{F}_{n}\right)  $ into following sets:
\begin{align*}
W_{C}  &  :=\left\{  \mu\in\mathsf{Fm}\left(  \mathbf{F}_{n}\right)  :\left|
\mathbf{F}_{n}/\mu\right|  =2\text{ and }C=\left\{  i=1,\ldots,n:\mathbf{x}%
_{i}/\mu\neq f\left(  \mathbf{x}_{i}\right)  /\mu\right\}  \right\} \\
&  \cup\left\{  \mu\in\mathsf{Fm}\left(  \mathbf{F}_{n}\right)  :\left|
\mathbf{F}_{n}/\mu\right|  >2\text{ and }C=\left\{  i=1,\ldots,n:\mathbf{x}%
_{i}/\mu=\ast_{\mu}\right\}  \right\}  \text{.}%
\end{align*}
for $C\subset\left\{  1,\ldots,n\right\}  $.

For $p\in T_{\mathcal{F}}(n)$ define
\[
N\left(  p\right)  :=\left\{  \mu\in\mathsf{Fm}\left(  \mathbf{F}_{n}\right)
:\mathbf{p}/\mu=\ast_{\mu}\right\}
\]
and
\[
k\left(  p\right)  :=\left|  \left\{  C:W_{C}\cap N\left(  p\right)
\neq\emptyset\right\}  \right|  \text{.}%
\]

Consider the family
\[
\mathcal{P}:=\left\{  p\in T_{\mathcal{F}}(n):\mathcal{V}\models p\left(
x_{1}yy,\ldots,x_{n}yy\right)  \approx p\left(  x_{1},\ldots,x_{n}\right)
yy\text{ and }f\left(  \mathbf{p}\right)  =\mathbf{1}\right\}
\]
Since $\mathbf{t}\in\mathcal{P}$, it is enough to prove that $p$ is
$\mathcal{V}$-projective for every $p\in\mathcal{P}$. We proceed by induction
on $k(p)$.

I. $k(p)=0$. Then $N\left(  p\right)  =\emptyset$, as $\bigcup\left\{
W_{C}:C\subset\left\{  1,\ldots,n\right\}  \right\}  =\mathsf{Fm}\left(
\mathbf{F}_{n}\right)  $, and so, by Proposition \ref{monsubirr}.2,
$\mathbf{p}=\mathbf{1}$.

II. $k(p)\geq1$. Let $C\subset\left\{  1,\ldots,n\right\}  $ be such that
$W_{C}\cap N\left(  p\right)  \neq\emptyset$. Define a substitution
$g_{C}:\left\{  x_{1},\ldots,x_{n}\right\}  \rightarrow T_{\mathcal{F}}(n)$
by
\[
g_{C}\left(  x_{i}\right)  :=\left\{
\begin{tabular}
[c]{ll}%
$x_{i}p$ & for $i\in C$\\
$x_{i}pp$ & for $i\notin C$%
\end{tabular}
\right.  \text{ .}%
\]
We show that $g_{C}\left(  p\right)  \in\mathcal{P}$ and $k\left(
g_{C}\left(  p\right)  \right)  <k(p)$.

Note that for the endomorphism $\mathbf{g}_{C}:\mathbf{F}_{n}\rightarrow
\mathbf{F}_{n}$ corresponding to $g_{C}$ we have $\mathbf{g}_{C}\left(
\mathbf{x}_{i}\right)  \equiv_{\left[  \mathbf{p}\right)  }\mathbf{x}_{i}$ for
$i=1,\ldots,n$, and so $\mathbf{g}_{C}\left(  \mathbf{p}\right)
\equiv_{\left[  \mathbf{p}\right)  }\mathbf{p}$. Hence $\mathbf{g}_{C}\left(
\mathbf{p}\right)  \in\left[  \mathbf{p}\right)  \subset\ker f$. To show that
$g_{C}\left(  p\right)  $ fulfills the identity from (1) take $\mathbf{A}%
\in\mathcal{V}$ and $a_{1},\ldots,a_{n},b\in A$. Put $\overline{p}%
:=p^{\mathbf{A}}\left(  a_{1},\ldots,a_{n}\right)  $. We know that
$p^{\mathbf{A}}\left(  a_{1}bb,\ldots,a_{n}bb\right)  =\overline{p}bb$. Then
$\left(  g_{C}\left(  p\right)  \right)  ^{\mathbf{A}}\left(  a_{1}%
bb,\ldots,a_{n}bb\right)  =p^{\mathbf{A}}\left(  c_{1},\ldots,c_{n}\right)  $,
where $c_{i}:=a_{i}bb\left(  \overline{p}bb\right)  =a_{i}\overline{p}bb$ for
$i\in C$ and $c_{i}:=a_{i}bb\left(  \overline{p}bb\right)  \left(
\overline{p}bb\right)  =a_{i}\overline{p}\overline{p}bb$ for $i\notin C$.
Hence\linebreak \ $\left(  g_{C}\left(  p\right)  \right)  ^{\mathbf{A}%
}\left(  a_{1}bb,\ldots,a_{n}bb\right)  =\left(  g_{C}\left(  p\right)
\right)  ^{\mathbf{A}}\left(  a_{1},\ldots,a_{n}\right)  bb$. This proves that
$g_{C}\left(  p\right)  \in\mathcal{P}$.

To prove that $k\left(  g_{C}\left(  p\right)  \right)  <k(p)$ we start from
the observation that \linebreak $\mathbf{g}_{C}\left(  \mathbf{p}\right)
\mathbf{pp}=\mathbf{g}_{C}\left(  \mathbf{p}\right)  $. From the definition of
$g_{C}$ we immediately obtain that $\mathbf{g}_{C}\left(  \mathbf{x}%
_{i}\right)  \mathbf{pp}=\mathbf{g}_{C}\left(  \mathbf{x}_{i}\right)  $ for
every $i=1,\ldots,n$. Hence and from the fact that $p\in\mathcal{P}$ we have
\begin{align*}
\mathbf{g}_{C}\left(  \mathbf{p}\right)  \mathbf{pp}  &  =\left(
\mathbf{g}_{C}\left(  p^{\mathbf{F}_{n}}\left(  \mathbf{x}_{1},\ldots
,\mathbf{x}_{n}\right)  \right)  \right)  \mathbf{pp}\\
&  =\left(  p^{\mathbf{F}_{n}}\left(  \mathbf{g}_{C}\left(  \mathbf{x}%
_{1}\right)  ,\ldots,\mathbf{g}_{C}\left(  \mathbf{x}_{n}\right)  \right)
\right)  \mathbf{pp}\\
&  =p^{\mathbf{F}_{n}}\left(  \mathbf{g}_{C}\left(  \mathbf{x}_{1}\right)
\mathbf{pp},\ldots,\mathbf{g}_{C}\left(  \mathbf{x}_{n}\right)  \mathbf{pp}%
\right) \\
&  =p^{\mathbf{F}_{n}}\left(  \mathbf{g}_{C}\left(  \mathbf{x}_{1}\right)
,\ldots,\mathbf{g}_{C}\left(  \mathbf{x}_{n}\right)  \right) \\
&  =\mathbf{g}_{C}\left(  \mathbf{p}\right)  \text{ .}%
\end{align*}

Now we show that $N\left(  g_{C}\left(  p\right)  \right)  \subset N(p)$. Let
$\mu\in N\left(  g_{C}\left(  p\right)  \right)  $. Then $\mathbf{g}%
_{C}\left(  \mathbf{p}\right)  /\mu=\ast_{\mu}$ and, since $\mathbf{g}%
_{C}\left(  \mathbf{p}\right)  \in\left[  \mathbf{p}\right)  $, we get
$\mathbf{p}\notin\mu$. On the other hand
\[
\left(  \mathbf{g}_{C}\left(  \mathbf{p}\right)  /\mu\right)  \left(
\mathbf{p}/\mu\right)  \left(  \mathbf{p}/\mu\right)  =\left(  \mathbf{g}%
_{C}\left(  \mathbf{p}\right)  \mathbf{pp}\right)  /\mu=\mathbf{g}_{C}\left(
\mathbf{p}\right)  /\mu=\ast_{\mu}\text{,}%
\]
and from Proposition \ref{monolith}.1 we get $\mathbf{p}/\mu=\ast_{\mu}$, as
required. Thus $k\left(  g_{C}\left(  p\right)  \right)  \leq k(p)$.

To prove that this inequality is sharp, it suffices to show that $W_{C}\cap
N\left(  g_{C}\left(  p\right)  \right)  =\emptyset$. To get a contradiction,
suppose that there exists $\mu\in W_{C}$ such that $\mathbf{g}_{C}\left(
\mathbf{p}\right)  /\mu=\ast_{\mu}$, and so, as we have just proved,
$\mathbf{p}/\mu=\ast_{\mu}$. We consider two cases: $\left|  \mathbf{F}%
_{n}/\mu\right|  =2$ and $\left|  \mathbf{F}_{n}/\mu\right|  >2$.

In the first case $C=\left\{  i=1,\ldots,n:\mathbf{x}_{i}/\mu\neq f\left(
\mathbf{x}_{i}\right)  /\mu\right\}  $. Then $\mathbf{g}_{C}\left(
\mathbf{x}_{i}\right)  /\mu=f\left(  \mathbf{x}_{i}\right)  /\mu$ for
$i=1,\ldots,n$, since it is easily seen that $\mathbf{g}_{C}\left(
\mathbf{x}_{i}\right)  /\mu$ $=\left(  \mathbf{x}_{i}/\mu\right)  \left(
\mathbf{p}/\mu\right)  =f\left(  \mathbf{x}_{i}\right)  /\mu$ for $i\in C$ and
$\mathbf{g}_{C}\left(  \mathbf{x}_{i}\right)  /\mu$ $=\left(  \mathbf{x}%
_{i}/\mu\right)  \left(  \mathbf{p}/\mu\right)  \left(  \mathbf{p}/\mu\right)
=\mathbf{x}_{i}/\mu=f\left(  \mathbf{x}_{i}\right)  /\mu$ for $i\notin C$.
Hence we get
\begin{align*}
\mathbf{g}_{C}\left(  \mathbf{p}\right)  /\mu &  =p^{\mathbf{F}_{n}}\left(
\mathbf{g}_{C}\left(  \mathbf{x}_{1}\right)  ,\ldots,\mathbf{g}_{C}\left(
\mathbf{x}_{n}\right)  \right)  /\mu\\
&  =p^{\mathbf{F}_{n}/\mu}\left(  \mathbf{g}_{C}\left(  \mathbf{x}_{1}\right)
/\mu,\ldots,\mathbf{g}_{C}\left(  \mathbf{x}_{n}\right)  /\mu\right) \\
&  =p^{\mathbf{F}_{n}/\mu}\left(  f\left(  \mathbf{x}_{1}\right)  /\mu
,\ldots,f\left(  \mathbf{x}_{1}\right)  /\mu\right) \\
&  =p^{\mathbf{F}_{n}}\left(  f\left(  \mathbf{x}_{1}\right)  ,\ldots,f\left(
\mathbf{x}_{n}\right)  \right)  /\mu\\
&  =f\left(  \mathbf{p}\right)  /\mu\\
&  =\mathbf{1}/\mu\text{, }%
\end{align*}
a contradiction.

In the second case $C=\left\{  i=1,\ldots,n:\mathbf{x}_{i}/\mu=\ast_{\mu
}\right\}  $. Then $\mathbf{g}_{C}\left(  \mathbf{x}_{i}\right)  /\mu\neq
\ast_{\mu}$ for $i=1,\ldots,n$, because $\mathbf{g}_{C}\left(  \mathbf{x}%
_{i}\right)  /\mu=\left(  \mathbf{x}_{i}/\mu\right)  \left(  \mathbf{p}%
/\mu\right)  =\mathbf{1}/\mu\neq\ast_{\mu}$ for $i\in C$ and, by Proposition
\ref{monolith}.1, $\mathbf{g}_{C}\left(  \mathbf{x}_{i}\right)  /\mu=\left(
\mathbf{x}_{i}/\mu\right)  \left(  \mathbf{p}/\mu\right)  \left(
\mathbf{p}/\mu\right)  =\mathbf{x}_{i}/\mu\neq\ast_{\mu}$ for $i\notin C$. On
the other hand, since $\left|  \mathbf{F}_{n}/\mu\right|  >2$ and
\[
\ast_{\mu}=\mathbf{g}_{C}\left(  \mathbf{p}\right)  /\mu=p^{\mathbf{F}_{n}%
/\mu}\left(  \mathbf{g}_{C}\left(  \mathbf{x}_{1}\right)  /\mu,\ldots
,\mathbf{g}_{C}\left(  \mathbf{x}_{n}\right)  /\mu\right)  \text{,}%
\]
we deduce from (2) that there is $i=1,\ldots,n$ such that $\mathbf{g}%
_{C}\left(  \mathbf{x}_{i}\right)  /\mu=\ast_{\mu}$, a contradiction.

Hence $k\left(  g_{C}\left(  p\right)  \right)  <k(p)$. Applying induction
assumption to $g_{C}\left(  p\right)  $ we get that $g_{C}\left(  p\right)  $
is $\mathcal{V}$-projective. So we find $\sigma\in\operatorname*{Hom}\left(
\mathbf{F}_{n},\mathbf{F}_{n}\right)  $ such that $\sigma\left(
\mathbf{g}_{C}\left(  \mathbf{p}\right)  \right)  =\mathbf{1}$, $\sigma\left(
\mathbf{x}_{i}\right)  \equiv_{\left[  \mathbf{g}_{C}\left(  \mathbf{p}%
\right)  \right)  }\mathbf{x}_{i} $, and so $\sigma\left(  \mathbf{g}%
_{C}\left(  \mathbf{x}_{i}\right)  \right)  \equiv_{\left[  \mathbf{g}%
_{C}\left(  \mathbf{p}\right)  \right)  }\mathbf{g}_{C}\left(  \mathbf{x}%
_{i}\right)  $ for every $i=1,\ldots,n$. Put $\tau:=\sigma\circ\mathbf{g}_{C}%
$. Then $\tau\left(  \mathbf{p}\right)  =\mathbf{1}$ and $\tau\left(
\mathbf{x}_{i}\right)  =\sigma\left(  \mathbf{g}_{C}\left(  \mathbf{x}%
_{i}\right)  \right)  \equiv_{\left[  \mathbf{g}_{C}\left(  \mathbf{p}\right)
\right)  }\mathbf{g}_{C}\left(  \mathbf{x}_{i}\right)  \equiv_{\left[
\mathbf{p}\right)  }\mathbf{x}_{i}$ for $i=1,\ldots,n$. Since $\left[
\mathbf{g}_{C}\left(  \mathbf{p}\right)  \right)  \subset\left[
\mathbf{p}\right)  $, we get $\tau\left(  \mathbf{x}_{i}\right)
\equiv_{\left[  \mathbf{p}\right)  }\mathbf{x}_{i}$ for $i=1,\ldots,n$, which
means that $\tau$ is a $\mathcal{V}$-projective unifier for $p$, as desired.
\end{proof}

\begin{remark}
In general, the reverse implication $(3)\nRightarrow(1)\Leftrightarrow(2)$
does not hold. Let $\mathcal{H}$ be the variety of Heyting algebras. Then
$t\left(  x_{1},x_{2}\right)  :=x_{2}\left(  x_{1}\vee\left(  x_{1}0\right)
\right)  $ is $\mathcal{H}$-projective (consider $\tau\in\operatorname*{Hom}%
\left(  \mathbf{F}_{2},\mathbf{F}_{2}\right)  $ given by $\tau\left(
\mathbf{x}_{1}\right)  :=\mathbf{x}_{1}$, $\tau\left(  \mathbf{x}_{2}\right)
:=\mathbf{x}_{1}\vee\left(  \mathbf{x}_{1}0\right)  $), but does not fulfil
$(1)$. To show this, take the subdirectly irreducible $5$-element nonlinear
Heyting algebra $\mathbf{H}:=\mathbf{2}^{\mathbf{2}}\oplus\mathbf{1}$ with the
universe $H=\left\{  0,a,a0,a\vee\left(  a0\right)  ,1\right\}  $. Then
$t^{\mathbf{H}}\left(  a00,100\right)  =t^{\mathbf{H}}\left(  a,1\right)
=a\vee\left(  a0\right)  \neq1$, whereas $t^{\mathbf{H}}\left(  a,1\right)
00=1$.
\end{remark}

It follows from Proposition \ref{matching}\ and Theorem \ref{123} that for a
congruence permutable Fregean variety $\mathcal{V}$ the sufficient condition
for having projective unifiers is that for every algebra $\mathbf{A}$ in
$\mathcal{V}$ and $a\in A$, the unary polynomial $\chi_{a}$ is an
endomorphism. The next theorem shows that this condition is also necessary.
This provides a full characterization of congruence permutable Fregean
varieties with projective unifiers. In the proof of the theorem we will use
the following computational lemma.

\begin{lemma}
\label{star}Let $\mathbf{A}$ be a subdirectly irreducible algebra in a Fregean
variety $\mathcal{V}$, $n\in\mathbb{N}$, $a_{1},\ldots,a_{n}\in A$, and let
$t\in T_{\mathcal{F}}(n)$ be a $\mathcal{V}$-projective term such that
$t^{\mathbf{A}}\left(  a_{1},\ldots,a_{n}\right)  =\ast$. Then $a_{i}=\ast$ or
$a_{i}=1$ for some $i\in\left\{  1,\ldots,n\right\}  $.
\end{lemma}

\begin{proof}
Suppose to the contrary that $a_{1},\ldots,a_{n}\in A\backslash\left\{
\ast,1\right\}  $. Since $t$ is projective there exists an endomorphism
$\tau:\mathbf{F}_{n}\rightarrow\mathbf{F}_{n}$ such that $\tau\left(
\mathbf{x}_{i}\right)  \equiv_{\left[  \mathbf{t}\right)  }\mathbf{x}_{i}$ and
$\tau\left(  \mathbf{t}\right)  =1$. Put $B:=\operatorname*{Sg}^{\mathbf{A}%
}\left(  a_{1},\ldots,a_{n}\right)  $ for the subalgebra of $\mathbf{A}$
generated by $\left\{  a_{1},\ldots,a_{n}\right\}  $. As $\ast\in B$, we get
$\Theta_{\mathbf{B}}\left(  \ast,1\right)  \subset\Theta_{\mathbf{A}}\left(
\ast,1\right)  |_{B}$. On the other hand, by Proposition \ref{FregSubIrr}.1,
we get $\Theta_{\mathbf{A}}\left(  \ast,1\right)  =\operatorname*{id}_{A}%
\cup\left\{  \left(  \ast,1\right)  ,\left(  1,\ast\right)  \right\}  $. Hence
$\Theta_{\mathbf{B}}\left(  \ast,1\right)  =\operatorname*{id}_{B}\cup\left\{
\left(  \ast,1\right)  ,\left(  1,\ast\right)  \right\}  $. Let us consider
the only epimorphism $f:\mathbf{F}_{n}\rightarrow\mathbf{B}$ such that
$f\left(  x_{i}\right)  =a_{i}$ for $i=1,\ldots,n$. Since $f\left(
\mathbf{t}\right)  =\ast$, we get $\left(  \mathbf{t},\mathbf{1}\right)  \in
f^{-1}\left(  \Theta_{\mathbf{B}}\left(  \ast,1\right)  \right)
\in\mathsf{Con}\left(  \mathbf{F}_{n}\right)  $. Hence $\Theta_{\mathbf{F}%
_{n}}\left(  \mathbf{t},\mathbf{1}\right)  \subset f^{-1}\left(
\Theta_{\mathbf{B}}\left(  \ast,1\right)  \right)  $. Let $i=1,\ldots,n$. Then
$\left(  \mathbf{x}_{i},\tau\left(  \mathbf{x}_{i}\right)  \right)  \in
f^{-1}\left(  \Theta_{\mathbf{B}}\left(  \ast,1\right)  \right)  $, and
therefore $\left(  a_{i},f\left(  \tau\left(  \mathbf{x}_{i}\right)  \right)
\right)  =\left(  f\left(  \mathbf{x}_{i}\right)  ,f\left(  \tau\left(
\mathbf{x}_{i}\right)  \right)  \right)  \in\Theta_{\mathbf{B}}\left(
\ast,1\right)  $. Hence, and since $a_{i}\notin\left\{  \ast,1\right\}  $, we
get $a_{i}=f\left(  \tau\left(  \mathbf{x}_{i}\right)  \right)  $. Finally, we
have $1=f\left(  \tau\left(  \mathbf{t}\right)  \right)  =f\left(
t^{\mathbf{F}_{n}}\left(  \tau\left(  \mathbf{x}_{1}\right)  ,\ldots
,\tau\left(  \mathbf{x}_{n}\right)  \right)  \right)  =t^{\mathbf{A}}\left(
f\left(  \tau\left(  \mathbf{x}_{1}\right)  \right)  ,\ldots,f\left(
\tau\left(  \mathbf{x}_{n}\right)  \right)  \right)  =t^{\mathbf{A}}\left(
a_{1},\ldots,a_{n}\right)  =\ast$, a contradiction.
\end{proof}

\begin{theorem}
\label{1234}Let $\mathcal{V}$ be a congruence permutable Fregean variety with
signature $\mathcal{F}$. Then the following conditions are equivalent:

\begin{enumerate}
\item [(1)]for every $f\in\mathcal{F}$%
\[
\mathcal{V}\models f\left(  x_{1}yy,\ldots,x_{k}yy\right)  \approx f\left(
x_{1},\ldots,x_{k}\right)  yy
\]
where $k$ is the arity of $f$;

\item[(2)] for every subdirectly irreducible $\mathbf{A}$ in $\mathcal{V}$, if
$\left|  A\right|  >2$, then $A\backslash\left\{  \ast\right\}  $ is a
subuniverse of $A$;

\item[(3)] $\mathcal{V}$ has projective unifiers.
\end{enumerate}
\end{theorem}

\begin{proof}
Clearly, it follows from Theorem \ref{123} and Proposition \ref{matching} that
it suffices to show that $\left(  3\right)  \Rightarrow\left(  2\right)  $. On
the contrary assume that there is a subdirectly irreducible $\mathbf{A}%
\in\mathcal{V}$ such that $\left|  A\right|  >2$ and there exist $t\in
T_{\mathcal{F}}(n)$, $n\in\mathbb{N}$, $a_{1},\ldots,a_{n}\in A\backslash
\left\{  \ast\right\}  $ with the property $t^{\mathbf{A}}\left(  a_{1}%
,\ldots,a_{n}\right)  =\ast$. Without loss of generality we can assume that
$\left\{  i:a_{i}\neq1\right\}  =\left\{  1,\ldots,m\right\}  $, where $0\leq
m\leq n$. Put $s\in T_{\mathcal{F}}(m+1)$ by
\[
s\left(  x_{1},\ldots,x_{m+1}\right)  :=t\left(  x_{1},\ldots,x_{m}%
,1,\ldots,1\right)  x_{m+1}x_{m+1}t\left(  x_{1},\ldots,x_{m},1,\ldots
,1\right)  \text{.}%
\]
Note that $s$ is unifiable, since for the substitution $\sigma\in
\operatorname*{Hom}\left(  \mathbf{F}_{m+1},\mathbf{F}_{m+1}\right)  $ given
by $\sigma\left(  \mathbf{x}_{i}\right)  :=\mathbf{x}_{i}$ for $i=1,\ldots,m$
and $\sigma\left(  \mathbf{x}_{m+1}\right)  :=\mathbf{1}$, we have
$\sigma\left(  \mathbf{s}\right)  =1$. Hence, by $(3)$, $s$ is projective. Let
$b\in A\backslash\left\{  \ast.1\right\}  $. Then
\begin{align*}
s^{\mathbf{A}}\left(  a_{1},\ldots,a_{m},b\right)   &  =t^{\mathbf{A}}\left(
a_{1},\ldots,a_{m},1,\ldots,1\right)  bbt^{\mathbf{A}}\left(  a_{1}%
,\ldots,a_{m},1,\ldots,1\right) \\
&  =t^{\mathbf{A}}\left(  a_{1},\ldots,a_{n}\right)  bbt^{\mathbf{A}}\left(
a_{1},\ldots,a_{n}\right) \\
&  =\ast bb\ast=\ast\text{.}%
\end{align*}
On the other hand, it follows from Lemma \ref{star} that $\left\{
a_{1},\ldots,a_{m},b\right\}  \cap\left\{  \ast.1\right\}  \neq\emptyset$, a contradiction.
\end{proof}

\begin{corollary}
For each congruence permutable Fregean variety $\mathcal{V}$ there exists the
largest subvariety that has projective unifiers. Namely, it is enough to put
\[
\mathcal{W}:=\left\{
\begin{array}
[c]{c}%
A\in\mathcal{V}:A\models f\left(  x_{1}yy,\ldots,x_{k}yy\right)  \approx
f\left(  x_{1},\ldots,x_{k}\right)  yy\text{ \ \ }\\
\text{for every }f\in\mathcal{F}\text{, where }k\text{ is the arity of }f
\end{array}
\right\}  \text{.}%
\]
\end{corollary}

\begin{remark}
In particular, for the variety of Heyting algebras $\mathcal{H}$, the largest
subvariety that has projective unifiers is
\[
\left\{  A\in\mathcal{H}:A\models x_{1}yy\vee x_{2}yy\approx\left(  x_{1}\vee
x_{2}\right)  yy\right\}  =\mathcal{LC}\text{,}%
\]
where $\mathcal{LC}$ denotes the variety of linear Heyting algebras (or
G\"{o}del-Dummett algebras) usually axiomatized by $\left(  x\rightarrow
y\right)  \vee\left(  y\rightarrow x\right)  \approx1$ (This fact was obtained
by straightforward calculation in $\mathcal{H}$ by Wro\'{n}ski, see
\cite{Wro08}).
\end{remark}

\begin{corollary}
From Theorem \ref{1234}, using condition (1) or (2), we can immediately deduce
some well known results. Namely, we get that

\begin{itemize}
\item  Boolean algebras ($\mathbf{CPC}$) (\cite{ButSim87},\cite{MarNip88},
\cite{MarNip89});

\item  Brouwerian semilattices ($\mathbf{IPC}$, $\rightarrow$, $\wedge$) (
\cite{Koh81}, see also \cite{NemWha71}, \cite{Ghi97});

\item  equivalential algebras ($\mathbf{IPC}$, $\leftrightarrow$)
(\cite{Wro05});

\item  Brouwerian semilattices with $0$ (bounded) ($\mathbf{IPC}$,
$\rightarrow$, $\wedge$, $\lnot$) (\cite{Wro05})

have projective unifiers.
\end{itemize}

The same is true for some other classes of algebras related to logic, e.g. for
equivalential algebras with $0$, i.e. the variety of algebras $\left(
A,\leftrightarrow,1,0\right)  $ denoted by $\mathcal{E}_{0}$, and being the
algebraic counterpart of the ($\leftrightarrow$, $\lnot$) fragment of
$\mathbf{IPC}$, in such a sense that for every $n\in\mathbb{N}$ and
$\varphi\in T_{\left(  \leftrightarrow,1,0\right)  }\left(  n\right)  $ we
have $\mathcal{E}_{0}\models\varphi\approx1$ iff $\vdash_{\mathbf{IPC}}%
\varphi$. Since we have here the identity $0xx\approx0$, the condition (1) of
Theorem \ref{1234} is fulfilled and so the variety has projective unifiers.
\end{corollary}

Under the additional assumption of finite signature, we can characterize
unifiable terms in the variety with projective unifiers.

\begin{proposition}
\label{unif}Let $\mathcal{V}$ be a congruence permutable Fregean variety with
finite signature $\mathcal{F}$ such that for every $f\in\mathcal{F}$%
\begin{equation}
\mathcal{V}\models f\left(  x_{1}yy,\ldots,x_{k}yy\right)  \approx f\left(
x_{1},\ldots,x_{k}\right)  yy \tag{*}%
\end{equation}
where $k$ is the arity of $f$, and let $t\in T_{\mathcal{F}}(n)$. Then the
following conditions are equivalent:

\begin{enumerate}
\item [(1)]$t$ is $\mathcal{V}$-unifiable;

\item[(2)] $t$ is $\mathcal{V}$-projective;

\item[(3)] for every $\mathbf{A}$ in $\mathcal{V}$, if $\left|  A\right|  =2$,
then $1^{\mathbf{A}}\in\operatorname{Im}t^{\mathbf{A}}$;

\item[(4)] $\left[  \mathbf{t}\right)  \cap\left\{  c^{\mathbf{F}_{n}}:c\in
T_{\mathcal{F}}\left(  0\right)  \right\}  =\left\{  \mathbf{1}\right\}  $.
\end{enumerate}
\end{proposition}

\begin{proof}
$(1)\Rightarrow(2)$. The implication follows from Theorem \ref{1234}.

$(2)\Rightarrow(3)$. Let $\mathbf{A}\in\mathcal{V}$ and let $\tau
\in\operatorname*{Hom}\left(  \mathbf{F}_{n},\mathbf{F}_{n}\right)  $ be a
$\mathcal{V}$-unifier\ for $t$. Take any homomorphism $\iota:\mathbf{F}%
_{n}\rightarrow\mathbf{A}$. Then $t^{\mathbf{A}}\left(  \iota\left(
\tau\left(  \mathbf{x}_{1}\right)  \right)  ,\ldots,\iota\left(  \tau\left(
\mathbf{x}_{n}\right)  \right)  \right)  =\iota\left(  t^{\mathbf{F}_{n}%
}\left(  \tau\left(  \mathbf{x}_{1}\right)  ,\ldots,\tau\left(  \mathbf{x}%
_{n}\right)  \right)  \right)  =\iota\left(  \tau\left(  \mathbf{t}\right)
\right)  =\iota\left(  \mathbf{1}\right)  =1^{\mathbf{A}}$.

$(3)\Rightarrow(4)$. Let $c\in T_{\mathcal{F}}\left(  0\right)  $ fulfill
$\mathbf{c:=}c^{\mathbf{F}_{n}}\in\left[  \mathbf{t}\right)  $. Suppose that
$\mathbf{c}\neq\mathbf{1}$. Then, by Proposition \ref{monsubirr}.2, there is
$\mu\in\mathsf{Fm}\left(  \mathbf{F}_{n}\right)  $ such that $\mathbf{c}%
/\mu=\ast_{\mu}$. Observe that $\left|  \mathbf{F}_{n}/\mu\right|  =2$, since
otherwise we would find $\mathbf{y}\notin\mu^{+}$ and, from Proposition
\ref{monolith}.1 and (*), we would get $\mathbf{c}/\mu=\left(  \mathbf{c}%
/\mu\right)  \left(  \mathbf{y}/\mu\right)  \left(  \mathbf{y}/\mu\right)
=\mathbf{1}/\mu$, a contradiction.

Put $\mathbf{A}:=\mathbf{F}_{n}/\mu$. It follows from (3) that there exist
$a_{1},\ldots,a_{n}\in A$ with $t^{\mathbf{A}}\left(  a_{1},\ldots
,a_{n}\right)  =1^{\mathbf{A}}$. Define a homomorphism $\pi:\mathbf{F}%
_{n}\rightarrow\mathbf{A}$ putting $\pi\left(  \mathbf{x}_{i}\right)  =a_{i}$
for $i=1,\ldots,n$. Then $\pi\left(  \mathbf{t}\right)  =\pi\left(
t^{\mathbf{F}_{n}}\left(  \mathbf{x}_{1},\ldots,\mathbf{x}_{n}\right)
\right)  =1^{\mathbf{A}}$ and so $\mathbf{c}\in\left[  \mathbf{t}\right)
\subset\ker\pi$. Thus $\mathbf{1}/\mu=1^{\mathbf{A}}=\pi\left(  \mathbf{c}%
\right)  =c^{\mathbf{A}}=\mathbf{c}/\mu=\ast_{\mu}$, contrary to our claim.

$(4)\Rightarrow(1)$. From Theorems \ref{DieGen} and \ref{1234} we deduce that
$\mathcal{V}$ is locally finite. Then $\Phi\left(  \mathbf{F}_{n}\right)  $ is
a finite modular lattice, which allows us to proceed by induction on $d\left(
\left[  \mathbf{t}\right)  \right)  $, where $d\left(  \left[  \mathbf{t}%
\right)  \right)  $ denotes the height of $\left[  \mathbf{t}\right)  $ in
$\Phi\left(  \mathbf{F}_{n}\right)  $. If $d\left(  \left[  \mathbf{t}\right)
\right)  =0$, $t$ is unifiable as $\mathbf{t}=\mathbf{1}$. Let $d\left(
\left[  \mathbf{t}\right)  \right)  >1$. Clearly, we can assume that
$t^{\mathbf{F}_{n}}\left(  \mathbf{1},\ldots,\mathbf{1}\right)  \neq
\mathbf{1}$, since otherwise $t$ is unifiable. From (4) we know that
$t^{\mathbf{F}_{n}}\left(  \mathbf{1},\ldots,\mathbf{1}\right)  \notin\left[
\mathbf{t}\right)  $. Thus we find a maximal $\mu\in\Phi\left(  \mathbf{F}%
_{n}\right)  $ such that $\left[  \mathbf{t}\right)  \subset\mu$ and
$t^{\mathbf{F}_{n}}\left(  \mathbf{1},\ldots,\mathbf{1}\right)  \notin\mu$. By
standard argument we get that $\mathbf{F}_{n}/\mu$ is subdirectly irreducible
and $t^{\mathbf{F}_{n}}\left(  \mathbf{1},\ldots,\mathbf{1}\right)  \in\mu
^{+}\backslash\mu$. From Proposition \ref{monsubirr}.1 we obtain $\ast_{\mu
}=t^{\mathbf{F}_{n}}\left(  \mathbf{1},\ldots,\mathbf{1}\right)
/\mu=t^{\mathbf{F}_{n}/\mu}\left(  \mathbf{1}/\mu,\ldots,\mathbf{1}%
/\mu\right)  $. From Theorem \ref{1234} we get immediately that $\left|
\mathbf{F}_{n}/\mu\right|  =2$. On the other hand, $\mathbf{1}/\mu
=\mathbf{t}/\mu=t^{\mathbf{F}_{n}}\left(  \mathbf{x}_{1},\ldots,\mathbf{x}%
_{n}\right)  /\mu=t^{\mathbf{F}_{n}/\mu}\left(  \mathbf{x}_{1}/\mu
,\ldots,\mathbf{x}_{n}/\mu\right)  $, and so $\left\{  i:\mathbf{x}_{i}%
/\mu=\ast_{\mu}\right\}  \neq\emptyset$. Without loss of generality we can
assume that $\left\{  i:\mathbf{x}_{i}/\mu=\ast_{\mu}\right\}  =\left\{
1,\ldots,m\right\}  $ for some $1\leq m\leq n$. Define $g:\left\{
x_{1},\ldots,x_{n}\right\}  \rightarrow T_{\mathcal{F}}(n)$ by
\[
g\left(  x_{i}\right)  :=\left\{
\begin{tabular}
[c]{ll}%
$x_{i}t$ & for $i=1,\ldots,m$\\
$x_{i}$ & for $i=m+1,\ldots,n$.
\end{tabular}
\right.  \text{ .}%
\]
Then for the endomorphism $\mathbf{g}:\mathbf{F}_{n}\rightarrow\mathbf{F}_{n}
$ corresponding to $g$ we have $\mathbf{x}_{i}\equiv_{\left[  \mathbf{t}%
\right)  }\mathbf{g}\left(  \mathbf{x}_{i}\right)  $ for $i=1,\ldots,n$. Hence
$\mathbf{t}\equiv_{\left[  \mathbf{t}\right)  }\mathbf{g}\left(
\mathbf{t}\right)  $, and finally $\mathbf{g}\left(  \mathbf{t}\right)
\in\left[  \mathbf{t}\right)  $. Now we show that $\mathbf{g}\left(
\mathbf{t}\right)  \neq\mathbf{t}$. In this aim, consider a homomorphism
$\psi:\mathbf{F}_{n}\rightarrow\mathbf{F}_{n}/\mu$ such that $\psi\left(
\mathbf{x}_{i}\right)  :=\mathbf{1}/\mu$, for $i=1,\ldots,n$. Then
$\psi\left(  \mathbf{t}\right)  =\psi\left(  t^{\mathbf{F}_{n}}\left(
\mathbf{x}_{1},\ldots,\mathbf{x}_{n}\right)  \right)  =t^{\mathbf{F}_{n}/\mu
}\left(  \mathbf{1}/\mu,\ldots,\mathbf{1}/\mu\right)  =\ast_{\mu}$. On the
other hand
\begin{align*}
\psi\left(  \mathbf{g}\left(  \mathbf{t}\right)  \right)   &  =\psi\left(
\mathbf{g}\left(  t^{\mathbf{F}_{n}}\left(  \mathbf{x}_{1},\ldots
,\mathbf{x}_{n}\right)  \right)  \right) \\
&  =\psi\left(  t^{\mathbf{F}_{n}}\left(  \mathbf{x}_{1}\mathbf{t}%
,\ldots,\mathbf{x}_{m}\mathbf{t},\mathbf{x}_{m+1}\mathbf{,}\ldots
,\mathbf{x}_{n}\right)  \right) \\
&  =t^{\mathbf{F}_{n}/\mu}\left(  \ast_{\mu},\ldots,\ast_{\mu},\mathbf{1}%
/\mu\mathbf{,}\ldots,\mathbf{1}/\mu\right) \\
&  =t^{\mathbf{F}_{n}/\mu}\left(  \mathbf{x}_{1}/\mu,\ldots,\mathbf{x}_{n}%
/\mu\right) \\
&  =t^{\mathbf{F}_{n}}\left(  \mathbf{x}_{1},\ldots,\mathbf{x}_{n}\right)
/\mu\\
&  =\mathbf{t}/\mu=1\text{ .}%
\end{align*}
Hence $\mathbf{g}\left(  \mathbf{t}\right)  \in\ker\psi$, whereas
$\mathbf{t}\notin\ker\psi$. Thus $\mathbf{t}<\mathbf{g}\left(  \mathbf{t}%
\right)  $, and so $d\left(  \left[  \mathbf{g}\left(  \mathbf{t}\right)
\right)  \right)  <d\left(  \left[  \mathbf{t}\right)  \right)  $. Moreover,
$\left[  \mathbf{g}\left(  \mathbf{t}\right)  \right)  \cap\left\{
c^{\mathbf{F}_{n}}:c\in T_{\mathcal{F}}\left(  0\right)  \right\}  =\left\{
\mathbf{1}\right\}  $, and from the induction hypothesis we know that
$g\left(  t\right)  $ is unifiable, i.e. there is $\sigma:\mathbf{F}%
_{n}\rightarrow\mathbf{F}_{m}$, $m\in\mathbb{N}$, a $\mathcal{V}$-unifier for
$g\left(  t\right)  $. Then $\sigma\circ\mathbf{g}$ is a $\mathcal{V}$-unifier
for $t$, which completes the proof.\smallskip
\end{proof}

As a corollary we get another sufficient and necessary condition for a
congruence permutable Fregean variety to have projective unifiers, which
guarantees that in such a variety every finitely generated algebra for which
any two constants are not glued together is projective.

\begin{corollary}
\label{fingen}Let $\mathcal{V}$ be a congruence permutable Fregean variety
with finite signature $\mathcal{F}$. Then $\mathcal{V}$ has projective
unifiers iff for every $n\in\mathbb{N}$ and $\varphi\in\Phi\left(
\mathbf{F}_{n}\right)  $%
\[
\left(  \varphi\cap\left\{  c^{\mathbf{F}_{n}}:c\in T_{\mathcal{F}%
}(0)\right\}  =\left\{  \mathbf{1}\right\}  \right)  \Rightarrow\left(
\mathbf{F}_{n}/\varphi\text{ is projective}\right)  \text{.}%
\]
\end{corollary}

\begin{proof}
$\Rightarrow)$ Let $\varphi\in\Phi\left(  \mathbf{F}_{n}\right)  $ satisfies
the condition $\varphi\cap\left\{  c^{\mathbf{F}_{n}}:c\in T_{\mathcal{F}%
}(0)\right\}  =\left\{  \mathbf{1}\right\}  $. It follows for Theorems
\ref{DieGen} and \ref{1234} that $\mathcal{V}$ is locally finite. Hence
$\varphi=\bigvee_{j=1}^{k}\left[  \mathbf{p}_{j}\right)  $ for some
$p_{1},\ldots,p_{k}\in T_{\mathcal{F}}(n)$. For every $\mathbf{s}\in\varphi$
we have\linebreak \ $\left[  \mathbf{s}\right)  \cap\left\{  c^{\mathbf{F}%
_{n}}:c\in T_{\mathcal{F}}(0)\right\}  =\left\{  \mathbf{1}\right\}  $, and
so, from Theorem \ref{1234}\ and Proposition \ref{unif}, we deduce that $s$ is
$\mathcal{V}$-unifiable. Now, from Lemma \ref{proj} we get the assertion.

$\Leftarrow)$ Let $t\in T_{\mathcal{F}}(n)$ be unifiable. Then we can find
$\sigma\in\operatorname*{Hom}\left(  \mathbf{F}_{n},\mathbf{F}_{n}\right)  $
such that $\left[  \mathbf{t}\right)  \subset\ker\sigma$. For $c^{\mathbf{F}%
_{n}}\in\left[  \mathbf{t}\right)  $ we have $c^{\mathbf{F}_{n}}=\sigma\left(
c^{\mathbf{F}_{n}}\right)  =\mathbf{1}$ for $c\in T_{\mathcal{F}}(0)$, and so
$c^{\mathbf{F}_{n}}=\mathbf{1}$. Thus $\mathbf{F}_{n}/\left[  \mathbf{t}%
\right)  $ is projective and the assertion follows from Proposition
\ref{projalg}.\smallskip
\end{proof}

It follows from Theorem \ref{1234} that as long as a congruence permutation
Fregean variety $\mathcal{V}$ has projective unifiers, every constant term $c
$ must be regular, that is $\mathcal{V}\models c\approx cyy$. In particular,
we can use our results to solve the \textsl{equational unification problem
with constants} in $\mathcal{V}$ such that $\mathcal{V}\models x\approx xyy$,
i.e. the equivalential reducts of algebras in $\mathcal{V}$ are Boolean
groups.\ It is for instance the case for the variety of Boolean algebras.

\section{Subtractive Fregean varieties}

Congruence permutability in Theorem \ref{1234} cannot be replaced by
subtractivity. To show this, consider the variety $\mathcal{HI}_{0}$ of
bounded Hilbert algebras, being the algebraic semantics of $(\rightarrow
,\lnot)$-reduct of $\mathbf{IPC}$, which is Fregean and subtractive but not
congruence permutable. Clearly, $\mathcal{HI}_{0}$ fulfills condition (2) from
Theorem \ref{1234}, but though the term $t\left(  x,y\right)  :=x\rightarrow
\left(  y\rightarrow0\right)  \in T_{\left(  \rightarrow,0\right)  }(2)$ is
unifiable in $\mathcal{HI}_{0}$, yet it has no mgus and, in consequence, it
has no projective unifiers, see \cite{Wro95}. Nevertheless, one can prove
weaker versions of Theorems \ref{123}, \ref{1234} and Corollary \ref{fingen}
valid for subtractive varieties. Let us start from the following proposition,
being an analogue of Theorems \ref{123}, that gives the sufficient condition
for \linebreak $\mathcal{V}$-projectivity.

\begin{proposition}
\label{subpro}If $\mathcal{V}$ is a locally finite subtractive Fregean variety
and \linebreak $t\in T_{\mathcal{F}}(n)$ fulfills the following two conditions:

\begin{itemize}
\item  for every subdirectly irreducible $\mathbf{A}$ in $\mathcal{V}$ with
$\left|  A\right|  >2$, and for all \linebreak $a_{1},\ldots,a_{n}\in A$%
\[
\left(  t^{\mathbf{A}}\left(  a_{1},\ldots,a_{n}\right)  =\ast\right)
\Rightarrow\left(  \exists i=1,\ldots,n:a_{i}=\ast\right)  \text{;}%
\]

\item $t$ is $\mathcal{V}$-unifiable with the special $\mathcal{V}$-unifier
$\tau\in\operatorname*{Hom}\left(  \mathbf{F}_{n},\mathbf{F}_{n}\right)  $
given by $\tau\left(  \mathbf{x}_{i}\right)  :=\mathbf{1}$ for $i=1,\ldots,n$
(i.e. $\mathcal{V}\models t\left(  1,\ldots,1\right)  \approx1$),
\end{itemize}

then $t$ is $\mathcal{V}$-projective.
\end{proposition}

\begin{proof}
Since $\mathcal{V}$ is locally finite, we can prove the assertion by induction
on $d\left(  \left[  \mathbf{t}\right)  \right)  $, where $d\left(  \left[
\mathbf{t}\right)  \right)  $ denotes the height of $\left[  \mathbf{t}%
\right)  $ in a finite modular lattice $\Phi\left(  \mathbf{F}_{n}\right)  $.
If $d\left(  \left[  \mathbf{t}\right)  \right)  =0$, then $\mathbf{t}%
=\mathbf{1} $, and so $t$ is $\mathcal{V}$-projective. Let $d\left(  \left[
\mathbf{t}\right)  \right)  =k\geq1$. Then $\mathbf{t}\neq\mathbf{1}$ and, by
Proposition \ref{monsubirr}.2, we can find $\mu\in\mathsf{Fm}\left(
\mathbf{F}_{n}\right)  $ such that $\mathbf{t}/\mu=\ast_{\mu}$. Put
$S:=\left\{  i=1,\ldots,n:\mathbf{x}_{i}/\mu=\ast_{\mu}\right\}  $. Define
$g:\left\{  x_{1},\ldots,x_{n}\right\}  \rightarrow T_{\mathcal{F}}(n)$ with
aid of the subtractive term $\mathsf{s}$ by
\[
g\left(  x_{i}\right)  =\left\{
\begin{tabular}
[c]{ll}%
$x_{i}$ & for $i\notin S$\\
$\mathsf{s}(t,x_{i})$ & for $i\in S$%
\end{tabular}
\right.  \text{ .}%
\]
Then $\mathsf{s}^{\mathbf{F}_{n}}(\mathbf{t},\mathbf{x}_{i})\equiv_{\left[
\mathbf{t}\right)  }\mathsf{s}^{\mathbf{F}_{n}}(\mathbf{1},\mathbf{x}%
_{i})=\mathbf{x}_{i}$, and so for the endomorphism $\mathbf{g}:\mathbf{F}%
_{n}\rightarrow\mathbf{F}_{n}$ corresponding to $g$ we have $\mathbf{g}\left(
\mathbf{x}_{i}\right)  \equiv_{\left[  \mathbf{t}\right)  }\mathbf{x}_{i} $
for $i=1,\ldots,n$. Hence $\mathbf{g}\left(  \mathbf{t}\right)  \equiv
_{\left[  \mathbf{t}\right)  }\mathbf{t}$, and, in consequence, $\mathbf{g}%
\left(  \mathbf{t}\right)  \in\left[  \mathbf{t}\right)  $ or, in other words,
$\mathbf{t}\leq\mathbf{g}\left(  \mathbf{t}\right)  $. We prove that
$\mathbf{t}<\mathbf{g}\left(  \mathbf{t}\right)  $.

Note that $\mathbf{g}\left(  \mathbf{x}_{i}\right)  /\mu=\mathsf{s}%
^{\mathbf{F}_{n}}(\mathbf{t},\mathbf{x}_{i})/\mu=\mathsf{s}^{\mathbf{F}%
_{n}/\mu}(\mathbf{t}/\mu,\mathbf{x}_{i}/\mu)=\mathsf{s}^{\mathbf{F}_{n}/\mu
}(\ast_{\mu},\ast_{\mu})=\mathbf{1}/\mu$ for $i\in S$. Hence $\mathbf{g}%
\left(  \mathbf{x}_{i}\right)  /\mu\neq\ast_{\mu}$ for $i=1,\ldots,n$. From
our assumptions one can easily deduce that for every subdirectly irreducible
$\mathbf{A}$ in $\mathcal{V}$, if \linebreak $a_{1},\ldots,a_{n}\in
A\backslash\left\{  \ast\right\}  $, then $t^{\mathbf{A}}\left(  a_{1}%
,\ldots,a_{n}\right)  \neq\ast$. Thus $\mathbf{g}\left(  \mathbf{t}\right)
/\mu=$\linebreak $\mathbf{t}^{\mathbf{F}_{n}/\mu}\left(  \mathbf{g}\left(
\mathbf{x}_{1}\right)  /\mu,\ldots,\mathbf{g}\left(  \mathbf{x}_{n}\right)
/\mu\right)  \neq\ast_{\mu}$, which means that $\mathbf{g}\left(
\mathbf{t}\right)  \neq\mathbf{t}$ since $\mathbf{t}/\mu=\ast_{\mu}$.
Therefore $d\left(  \left[  \mathbf{g}\left(  \mathbf{t}\right)  \right)
\right)  <d\left(  \left[  \mathbf{t}\right)  \right)  $. Applying the
induction hypothesis we know that $g\left(  t\right)  $ is $\mathcal{V}%
$-projective. Hence there exists an endomorphism $\tau:\mathbf{F}%
_{n}\rightarrow\mathbf{F}_{n}$ such that $\tau\left(  \mathbf{x}_{i}\right)
\equiv_{\left[  \mathbf{g}\left(  \mathbf{t}\right)  \right)  }\mathbf{x}_{i}$
for $i=1,\ldots,n$, and $\tau\left(  \mathbf{g}\left(  \mathbf{t}\right)
\right)  =1$. To complete the proof take an endomorphism $\tau\circ
\mathbf{g}:\mathbf{F}_{n}\rightarrow\mathbf{F}_{n}$. Then $\tau\left(
\mathbf{g}\left(  \mathbf{x}_{i}\right)  \right)  \equiv_{\left[
\mathbf{g}\left(  \mathbf{t}\right)  \right)  }\mathbf{g}\left(
\mathbf{x}_{i}\right)  \equiv_{\left[  \mathbf{t}\right)  }\mathbf{x}_{i} $,
which implies $\tau\left(  \mathbf{g}\left(  \mathbf{x}_{i}\right)  \right)
\equiv_{\left[  \mathbf{t}\right)  }\mathbf{x}_{i}$ for $i=1,\ldots,n$, and so
$\tau\circ\mathbf{g}$ is a $\mathcal{V}$-projective unifier for $t$, as required.\smallskip
\end{proof}

The final theorem presented below characterizes those subtractive Fregean
varieties with finite signature for which all terms (and hence all finite sets
of equations) are projective. As a corollary, we deduce that each such variety
is primitive, i.e. every its subquasivariety is a variety.

\begin{theorem}
\label{subthe}Let $\mathcal{V}$ be a subtractive Fregean variety with the
finite signature $\mathcal{F}$. Then the following conditions are equivalent:

\begin{enumerate}
\item [(1)]all terms in $T_{\mathcal{F}}(\omega)$ are $\mathcal{V}$-projective;

\item[(2)] for every subdirectly irreducible $\mathbf{A}$ in $\mathcal{V}$ the
subset $A\backslash\left\{  \ast\right\}  $ is a subuniverse of $A$;

\item[(3)] every finitely generated algebra in $\mathcal{V}$ is projective;

\item[(4)] every finite set of equations in $T_{\mathcal{F}}(\omega)$ is
$\mathcal{V}$-projective.
\end{enumerate}
\end{theorem}

\begin{proof}
$(1)\Rightarrow(2)$. To obtain the contradiction, assume that there is a
subdirectly irreducible $\mathbf{A}\in\mathcal{V}$ and $t\in T_{\mathcal{F}%
}(n)$, $n\in\mathbb{N}$, $a_{1},\ldots,a_{n}\in A\backslash\left\{
\ast\right\}  $ such that $t^{\mathbf{A}}\left(  a_{1},\ldots,a_{n}\right)
=\ast$. Without loss of generality we can assume that $\left\{  i:a_{i}%
\neq1\right\}  =\left\{  1,\ldots,k\right\}  $, where $0\leq k\leq n$. Put
$p\in T_{\mathcal{F}}(k)$ by $p\left(  x_{1},\ldots,x_{k}\right)  :=t\left(
x_{1},\ldots,x_{k},1,\ldots,1\right)  $. Then $p$ is $\mathcal{V}$-projective
and $p\left(  a_{1},\ldots,a_{k}\right)  =\ast$. It follows from Lemma
\ref{star} that there is $i=1,\ldots,k$ such that $a_{i}\in\left\{
\ast,1\right\}  $, a contradiction.

$(2)\Rightarrow(1)$. It follows from Theorems \ref{DieGen} that $\mathcal{V}$
is locally finite. Let $t\in T_{\mathcal{F}}(n)$. By Proposition \ref{subpro}
it is enough to show that $t^{\mathbf{F}_{n}}\left(  \mathbf{1},\ldots
,\mathbf{1}\right)  =\mathbf{1}$. On the contrary assume that $t^{\mathbf{F}%
_{n}}\left(  \mathbf{1},\ldots,\mathbf{1}\right)  \neq\mathbf{1}$. Then, using
Proposition \ref{monsubirr}.2, we find $\mu\in\mathsf{Fm}\left(
\mathbf{F}_{n}\right)  $ such that $t^{\mathbf{F}_{n}}\left(  \mathbf{1}%
,\ldots,\mathbf{1}\right)  /\mu=\ast_{\mu}$. Hence $t^{\mathbf{F}_{n}/\mu
}\left(  \mathbf{1}/\mu,\ldots,\mathbf{1}/\mu\right)  =t^{\mathbf{F}_{n}%
}\left(  \mathbf{1},\ldots,\mathbf{1}\right)  /\mu=\ast_{\mu}$, a contradiction.

$(1)$ and $(2)\Rightarrow(3)$. It is enough to show that $\mathbf{F}%
_{n}/\varphi$ is projective for every $n\in\mathbb{N}$ and $\varphi\in
\Phi\left(  \mathbf{F}_{n}\right)  $. It follows from (2) and Theorem
\ref{DieGen} that $\mathcal{V}$ is locally finite. Hence $\varphi
=\bigvee_{j=1}^{k}\left[  \mathbf{p}_{j}\right)  $ for some $p_{1}%
,\ldots,p_{k}\in T_{\mathcal{F}}(n)$. From (1) and Lemma \ref{proj} we obtain
the assertion.

$(3)\Rightarrow(4)$. Consider a set of equations $\left\{  s_{i}%
=t_{i}:i=1,\ldots,k\right\}  $, where $s_{i},t_{i}\in T_{\mathcal{F}}(n)$,
$i=1,\ldots,k$, for some $n\in\mathbb{N}$. Put $\varphi:=\bigvee_{j=1}%
^{k}\Theta\left(  \mathbf{s}_{j},\mathbf{t}_{j}\right)  $. As $\mathbf{F}%
_{n}/\varphi$ is finitely generated, it follows from (3) that $\mathbf{F}%
_{n}/\varphi$ is projective. Now, from Proposition \ref{projquotient} we get
that $\left\{  s_{i}=t_{i}:i=1,\ldots,k\right\}  $ is $\mathcal{V}%
$-projective, as desired.

$(4)\Rightarrow(1)$. The implication is obvious.\smallskip
\end{proof}

\begin{corollary}
If $\mathcal{V}$ is a subtractive Fregean variety with the finite signature
$\mathcal{F}$ that satisfies one of the equivalent conditions of Theorem
\ref{subthe}, then $\mathcal{V}$ is a primitive variety.
\end{corollary}

\begin{proof}
The proof is similar to one given in \cite[Corollary 4.7]{Slo96} for
$\mathcal{E}$, see also \cite{OlsRaf07}. As usual, by $H$, $S$ and $P_{U}$ we
denote the operations of the formation of homomorphic images, subalgebras and
ultraproducts, respectively. Let $\mathcal{Q}$ be a subvariety of
$\mathcal{V}$. It suffices to show that $H\left(  \mathcal{Q}\right)
\subset\mathcal{Q}$. Let $\mathbf{A}\in\mathcal{Q}$, $\mathbf{B}\in H\left(
\mathbf{A}\right)  $ and $\pi:\mathbf{A}\rightarrow\mathbf{B}$ be an
epimorphism. It is well known that $\mathbf{B}$ can be embedded into an
ultraproduct of its finitely generated subalgebras. On the other hand, it
follows from condition (3) of Theorem \ref{subthe} that each such subalgebra
$\mathbf{C}$ is projective, and so there exist a homomorphism $\iota
:\mathbf{C}\rightarrow\mathbf{A}$\textbf{\ }such that $\pi\circ\iota
=\operatorname*{id}_{C}$. Clearly, $\iota$ is a monomorphism and so
$\mathbf{C}\in S\left(  \mathbf{A}\right)  $. Since quasivarieties are closed
under the formation of subalgebras and ultraproducts, we get $\mathbf{B}\in
SP_{U}S\left(  \mathbf{A}\right)  \subset\mathcal{Q}$, as desired.\smallskip
\end{proof}

\end{document}